\documentclass{article}

\usepackage{amsthm,amssymb}
\catcode`\@=11 \@addtoreset{equation}{section}

\newtheorem{thm}{Theorem}[section]
\newtheorem{prop}{Proposition}[section]
\newtheorem{cor}{Corollary}[section]
\newtheorem{rem}{Remark}[section]
\newtheorem{lem}{Lemma}[section]
\newtheorem{defin}{Definition}[section]

\newcommand{\eps}{\varepsilon}
\newcommand{\loc}{{\mathrm{loc}}}

\newcommand{\rank}{{\mathrm{rank}\,}}

\newcommand{\const}{\mathrm{const}}
\newcommand{\R}{{\mathbb{R}}}
\newcommand{\cbf}{{\mathbf{c}}}
\newcommand{\Cbf}{{\mathbf{C}}}

\newcommand{\vbf}{{\mathbf{v}}}
\newcommand{\ubf}{{\mathbf{u}}}

\newcommand{\Xbf}{{\mathbf{X}}}

\newcommand{\zbf}{{\mathbf{z}}}

\newcommand{\C}{{\mathbb{C}}}
\newcommand{\Z}{{\mathbb{Z}}}

\newcommand{\calA}{\mathcal{A}}

\newcommand{\calB}{\mathcal{B}}

\newcommand{\calH}{\mathcal{H}}
\newcommand{\calL}{\mathcal{L}}
\newcommand{\calM}{\mathcal{M}}
\newcommand{\calN}{\mathcal{N}}
\newcommand{\calF}{\mathcal{F}}
\newcommand{\calP}{\mathcal{P}}

\newcommand{\calR}{\mathcal{R}}

\newcommand{\calT}{\mathcal{T}}
\newcommand{\calV}{\mathcal{V}}
\newcommand{\calW}{\mathcal{W}}
\newcommand{\calK}{\mathcal{K}}

\newcommand{\Crit}{{\mathrm{Crit}}}

\title{Shilnikov Lemma for a nondegenerate critical manifold of a
Hamiltonian system.
\author{Sergey Bolotin\thanks{Supported by the Programme
``Dynamical Systems and Control Theory'' of RAS and  RFBR grants \#12-01-00441 and \#13-01-12462.}
 \\
University of Wisconsin--Madison\\ and\\ Moscow Steklov Mathematical Institute \and
Piero Negrini\\ Department of Matematics
\\ Sapienza, University of Rome} }

\date{}

\begin{document}

\maketitle

\begin{abstract}
We prove an analog of Shilnikov Lemma   for a normally hyperbolic symplectic critical
manifold $M\subset H^{-1}(0)$ of a Hamiltonian system. Using this result,   trajectories
with small energy $H=\mu>0$ shadowing chains of
homoclinic orbits to $M$ are represented as extremals of a discrete variational problem,
and their existence is proved. This paper is motivated by
applications to the Poincar\'e second species
 solutions of the 3 body problem with 2 masses small of order $\mu$. As $\mu\to 0$,
double collisions of small bodies
correspond to a symplectic critical manifold of the regularized Hamiltonian system.
\end{abstract}

\section{Introduction}

Consider a smooth Hamiltonian system $(\calM,\omega,H)$ with phase space $ \calM$, symplectic form $\omega$ and Hamiltonian $H$.  Let $\vbf=\vbf_H$ be
the Hamiltonian vector field: $\omega(\vbf(x),\cdot )=-dH(x)$, and $\phi^t=\phi_H^t$ the flow of the system. Suppose that
$H$ has a nondegenerate normally hyperbolic symplectic critical $2m$-dimensional manifold $M\subset \Sigma_0=H^{-1}(0)$
with real eigenvalues.
Thus for any $z\in M$:
\begin{itemize}
\item $\rank d^2H(z)=2k=\dim\calM-2m$;
\item the restriction $\omega|_{T_zM}$   is nondegenerate;
\item the eigenvalues of the linearization of $\vbf$ at $z$ are all real.
\end{itemize}

Let  $D\phi^t(z)=e^{tA(z)}$ be the linearized flow. Denote by
$$
E_z=T_z^\perp M=\{\xi\in T_z \calM:\omega(\xi,\eta)=0\;\mbox{for all}\; \eta\in T_zM\}
$$
the symplectic complement to $T_z M $.
Since $M$ is symplectic, $T_z \calM =T_zM\oplus E_z$ and $\omega|_{E_z}$ is nondegenerate. Hence $E_z=E_z^+\oplus E_z^-$, where
$E_z^\pm$ are $k$-dimensional $A(z)$-invariant Lagrangian subspaces of $E_z$ corresponding to negative and positive  eigenvalues  respectively.
We write $\xi\in E_z$ as $\xi=(\xi_+,\xi_-)$, where $\xi_+\in E_z^+$ and $\xi_-\in E_z^-$.
Then the  linearized flow on $E_z$ is
\begin{equation}
\label{eq:linear}
D\phi^t(z)(\xi_+,\xi_-)=(e^{-tA_+(z)}\xi_+,e^{tA_-(z)}\xi_-),
\end{equation}
where the eigenvalues of $A_\pm(z)=\mp A(z)|_{E_z^\pm}$ are positive.
Thus $E_z^+$ is the stable subspace, and $E_z^-$ the unstable subspace.
The quadratic part of the Hamiltonian is
\begin{equation}
\label{eq:d2H} \frac12 d^2H(z)(\xi)=-\omega (\xi_-,A_+(z)\xi_+)=-\omega( A_-(z)\xi_-,\xi_+).
\end{equation}

The stable and unstable manifolds\footnote{In what follows $+$ corresponds to the stable manifold ($t\to +\infty$), and $-$
to the unstable manifold ($t\to -\infty$).}
$$
W^{\pm}(z)=\{x\in  \calM :\lim_{t\to\pm\infty}\phi^t(x)=z\}
$$
of an equilibrium $z\in M$ have dimension $k$  and $T_zW^\pm(z)=E_z^\pm$.  The stable and unstable manifolds
$$
W^\pm=W^\pm(M)=\cup_{z\in M}W^{\pm}(z)
$$
of  $M$ have dimension $k+2m$ and $T_zW^\pm=T_zM\oplus E_z^\pm$ for
any $z\in M$. It is well known (see e.g.\ \cite{DLS}) that  $W^\pm(z)$ are isotropic: $\omega|_{W^\pm(z)}=0$,
 and $W^\pm$ are coisotropic: for any $a\in W^\pm(z)$, we have $T_a^\perp W^\pm=T_aW^\pm(z)$.
Thus $W^\pm(z)$  form a smooth  isotropic foliation of $W^\pm$.
Define projections $\pi_\pm:W^\pm\to M$  by $\pi_\pm(x)=z$ if $x\in W^\pm(z)$:
$$
\pi_\pm(x)=\lim_{t\to\pm\infty} \phi^t(x).
$$

Since $M\subset \Sigma_0=H^{-1}(0)$, we have $W^\pm\subset\Sigma_0$. The intersection $\Gamma=(W^+\cap W^-)\setminus M$ consists of  orbits $\gamma:\R\to \calM$ homoclinic
to $M$, i.e.\ heteroclinic from $z_-=\gamma(-\infty)\in M$ to $z_+=\gamma(+\infty)\in M$.
Define a {\it scattering map} $\calF:\pi_-(\Gamma)\to \pi_+(\Gamma)$ setting $\calF(z_-)=z_+$ if there is an orbit  heteroclinic from $z_-$ to $z_+$, i.e.\   $W^-(z_-)\cap W^+(z_+)\ne \emptyset$.

\begin{rem}
Following  \cite{DLS}, we call $\calF$ the scattering map. However,
our case is different from  \cite{DLS} because the
manifold $M$ is critical. In particular, there is no straightforward
cross section for the flow near $M$.   The scattering map is also called the
homoclinic map. In the applications to Celestial
Mechanics \cite{Bol:ELL,Bol-Neg:DCDS}, we call $\calF$ the collision map.
\end{rem}

In general $\calF$ is multivalued.
To define a single valued smooth map, we need to consider local branches of $\calF$.
We call a heteroclinic  orbit $\gamma(t)=\phi^t(a)$, $\gamma(\pm\infty)=c_\pm\in M$, {\it transverse}
if the following conditions hold.

\begin{prop}\label{prop:coiso}
The following conditions are equivalent:
\begin{itemize}
\item $T_{a}W^-(c_-)\cap T_{a}W^+= \R \vbf(a)$,
\item $T_{a}W^+(c_+)\cap T_{a}W^-= \R \vbf(a)$,
\item The symplectic form $\omega$ defines a nondegenerate modulo $\R\vbf(a)$ bilinear form on
$T_{a}W^-(c_-)\times T_{a}W^+(c_+)$.
\item There exist Lagrangian submanifolds $\Lambda^\pm\subset M$ containing $c_\pm$ such that
the Lagrangian manifolds $W^\pm(\Lambda^\pm)=\cup_{z\in\Lambda^\pm}W^\pm(z)$  intersect transversely in $\Sigma_0$ along $\gamma$:
$$
T_{a}W^+(\Lambda^+)\cap T_{a}W^-(\Lambda^-)= \R \vbf(a).
$$
\end{itemize}
\end{prop}

These conditions imply that $a$ is a point of transverse intersection of $W^+$ and $W^-$, i.e.\
$T_{a}W^++T_{a}W^-=T_{a}\Sigma_0$. We skip an elementary proof of Proposition \ref{prop:coiso}.

If $\gamma$ is transverse, then $\calF$ has a well defined smooth branch $f:V^-\to V^+$,
where $V^\pm\subset M$ is a small neighborhood of $c_\pm$. Indeed, let $N\subset
W^+$  be a local section at $a$ such that $T_{a}N\oplus \R
\vbf(a)=T_{a}W^+$.  There exists a neighborhood $V^-\subset M$
of $c_-$ such that for any $z_-\in V^-$, the manifolds $W^-(z_-)$ and
$N$ intersect transversely in $\Sigma_0$ at a point $b$ close to $a$.  Set $z_+=f(z_-)=\pi_+(b)$.  Then
$\sigma(t)=\phi^t(b)$ is a heteroclinic orbit joining $z_-$ with $z_+$. The
map $f:V^-\to M$ is symplectic.

Indeed, let $(x_\pm,y_\pm)$ be local symplectic coordinates
in $V^\pm$ and $\alpha$ a 1-form in a neighborhood of $\overline{\gamma(\R)}$ such that
$d\alpha=\omega$ and $\alpha|_{V^\pm}=y_\pm\,dx_\pm$. Then by the first variation
formula \cite{Arnold}
$$
f(x_-,y_-)=(x_+,y_+)\quad \Rightarrow \quad
y_+\,dx_+-y_-\,dx_-=dG,\quad G(z_-)=\int_{\sigma}\alpha.
$$

We can choose symplectic coordinates $(x_\pm,y_\pm)$ in $V^\pm$ so that
$$
\Lambda^+=\{y_+=b_+\}=B(a_+)\times\{b_+\},\quad \Lambda^-=\{x_-=a_-\}=\{a_-\}\times B(b_-),
$$
where $c_\pm=(a_\pm,b_\pm)$ and $B$ is a small ball in $\R^m$.
Then for $(x_-,y_+)\in B(a_-)\times B(b_+)$, Lagrangian manifolds $W^-(\{x_-\}\times B(b_-))$
and $W^+(B(a_+)\times\{y_+\})$ intersect transversely in $\Sigma_0$ along a heteroclinic trajectory $\sigma(x_-,y_+)$
joining the points $(x_-,y_-)$ with $(x_+,y_+)$.
Decreasing the sets $V^\pm\subset M$ if necessary, we represent $f:V^-\to V^+$ by a generating function $S(x_-,y_+)=\langle y_+,x_+\rangle -G$ \cite{Arnold}:
\begin{equation}
\label{eq:genf}
f(x_-,y_-)=(x_+,y_+)\quad\Leftrightarrow\quad
 d S(x_-,y_+)=y_-\,dx_-+x_+\,dy_+.
\end{equation}

Introducing a local branch $f$ near any transverse heteroclinic orbit, we  represent the scattering map by a countable collection
$\calF$ of smooth  symplectic diffeomorphisms
$f:V^-\to V^+$ of open sets in $M$.
In general $\calF$ has infinitely many branches. For example, this is so in our
application to Celestial Mechanics \cite{Bol-Neg:DCDS}.  In fact $\calF$
being multivalued helps in constructing symbolic dynamics, see e.g.\ \cite{Bol:DCDS}.

An orbit of $\calF$ is a pair  of sequences $f_i:V_i^-\to V_i^+$ and
$z_i\in V_i=V_{i}^-\cap V_{i-1}^+$ such that $z_{i+1}=f_{i}(z_i)$.
It defines a chain  $\sigma=(\sigma_i)$ of transverse
heteroclinic orbits $\sigma_i$ connecting $z_i$ with
$z_{i+1}$.

\begin{rem}
The  scattering map may be viewed  as a single map -- the skew product of
 the maps $f\in\calF$ which is a (partly defined) map of $\calF^\Z\times M$.
This is needed to study  chaotic dynamics of $\calF$.
\end{rem}

Let $c_{i+1}=f_i(c_i)$ be a periodic orbit:
$f_{i+n}=f_i$, $c_{i+n}=c_i$.  Then $c_0$ is a fixed point of the composition
$F_n=f_{n-1}\circ\cdots\circ f_{0}$.
The periodic orbit is called {\it nondegenerate} if $z_0$ is a nondegenerate fixed point:
\begin{equation}
\label{eq:nondeg}
\det (DF_n(c_0)-I)\ne 0.
\end{equation}
Then the corresponding periodic heteroclinic chain $\sigma=(\sigma_i)$
will be called {\it nondegenerate}.

Let $z_i=(x_i,y_i)$ be symplectic coordinates in $V_i$ such that
$f_i$ is represented by a generating function
as in (\ref{eq:genf}):
\begin{equation}
\label{eq:genfj}
f_i(x_i,y_i)=(x_{i+1},y_{i+1})\quad\Leftrightarrow\quad
 d S_i(x_i,y_{i+1})=y_i\,dx_i+x_{i+1}\,dy_{i+1}.
\end{equation}
A periodic orbit of $\calF$ corresponds to a critical point $\cbf=(c_i)_{i=0}^{n-1}$  of the discrete action functional
\begin{equation}
\label{eq:calA}
\calA(\zbf)=\sum_{i=0}^{n-1}(S_i(x_i,y_{i+1})-\langle x_i,y_{i}\rangle),\qquad y_{n}=y_0.
\end{equation}
It is  well known (see \cite{Sympl}) that the periodic orbit is nondegenerate iff $\cbf$ is a nondegenerate critical point
of $\calA$.

To shadow a nondegenerate heteroclinic chain $\sigma$ by a trajectory of the
Hamiltonian system on $\Sigma_\mu=H^{-1}(\mu)$ with small $\mu\ne
0$, we need extra conditions which depend on the sign of $\mu$.

We assumed that the eigenvalues of equilibria in $M$ are real. There are
two main cases to consider:

\begin{itemize}
\item {\bf Generic real eigenvalues}: for any $z\in M$, eigenvalues of
    $A_\pm(z)$ satisfy
\begin{equation}
\label{eq:generic} 0<\lambda(z)=\lambda_1(z)<\lambda_2(z)\le\cdots \le \lambda_{k}(z).
\end{equation}
\item {\bf Equal semisimple eigenvalues}: for any $z\in M$,
\begin{equation}
\label{eq:Apm} A_\pm(z)=\lambda(z)I,\qquad \lambda(z)>0.
\end{equation}
\end{itemize}

The last case is highly nongeneric. However, it
appears in our main application \cite{Bol-Neg:DCDS} to Celestial Mechanics
 which is briefly discussed in the next section. For
this reason in this paper we  assume (\ref{eq:Apm}).
Generic real case is similar, but the details will be published elsewhere.
By (\ref{eq:linear})--(\ref{eq:d2H}) and (\ref{eq:Apm}),
\begin{eqnarray}
\label{eq:d2H2} d^2H(z)(\xi)&=&-2\lambda(z)\omega(\xi_-,\xi_+),\\
D\phi^t(z)(\xi)&=&(e^{-\lambda(z)t}\xi_+, e^{\lambda(z)t}\xi_-).
\end{eqnarray}

Since the flow on $W^\pm(z)$ is a node,
for any $a\in W^\pm(z)$ there exist tangent vectors
\begin{equation}
\label{eq:dir} \vbf_\pm(a)=\mp\lambda(z)\lim_{t\to\pm\infty}e^{\pm
t\lambda(z)}\vbf(\phi^t(a))\in E_{z}^\pm.
\end{equation}
The map $\vbf_\pm:W^\pm(z)\to E_z^\pm$ is smooth and $\vbf_\pm(z)=0$, $D\vbf_\pm(z)=I_{E_z^\pm}$
(see Proposition \ref{prop:change}).

\begin{rem}
In the case (\ref{eq:generic}) of generic real eigenvalues, $\vbf_\pm(a)=0$ for $a$ in the strong stable (unstable) manifold of $z$.
Otherwise, $\vbf_\pm(a)$ is collinear to the eigenvector $\ubf_\pm(z)$ of $A_\pm(z)$ associated to the eigenvalue
$\lambda(z)$.
\end{rem}

For a heteroclinic orbit $\gamma(t)=\phi^t(a)$ with $\gamma(\pm\infty)=z_\pm\in M$,
let $\vbf_\pm(\gamma)=\vbf_\pm(a)\in E_{z_\pm}^\pm$ be the vectors (\ref{eq:dir}).
They depend on the choice of the initial
point $a$ on $\gamma$, but the directions are well defined.

If $\sigma=(\sigma_i)$ is a heteroclinic chain, so that
$\sigma_{i-1}(+\infty)=\sigma_i(-\infty)=c_i\in M$, we set
\begin{equation}
\label{eq:change+} a_i(\sigma)=\omega(\vbf_i^+(\sigma),\vbf_i^-(\sigma)),\qquad \vbf_i^+(\sigma)=\vbf_+(\sigma_{i-1}),\quad \vbf_i^-(\sigma)=\vbf_-(\sigma_{i}).
\end{equation}

\begin{defin}\label{def:pos}
We call a heteroclinic  chain positive (negative) if $a_i(\sigma)>0$ $(a_i(\sigma)<0)$ for all $i$.
\end{defin}

\begin{rem}
This definition makes sense also for generic real eigenvalues. Then $v_i^\pm(\sigma)=k_i^\pm \ubf_\pm(c_i)$.
If we choose the eigenvectors $\ubf_\pm$ so that $\omega(\ubf_+,\ubf_-)>0$,
then the positivity condition means $k_i^-k_i^+>0$ for all $i$.
\end{rem}

Geometrically  the chain $\sigma$ is a piece wise smooth curve
$$
C=\cup_{i}\overline{\sigma_i(\R)}
$$
with ``reflections" from $M$  at the points $c_i$.
Then $a_i(\sigma)$ measures symplectic angles at these reflections.

Positive heteroclinic chains can be shadowed by orbits with small positive energy,
and negative chains with small negative energy. It is not possible  to shadow chains of mixed type.

\begin{thm}\label{thm:period}
Let   $\sigma$ be a positive nondegenerate periodic  heteroclinic
chain. Then there is $\mu_0>0$
such that for any $\mu\in (0,\mu_0]$:
\begin{itemize}
\item
There exists a periodic orbit
$\gamma_\mu$ on $\Sigma_\mu=H^{-1}(\mu)$, smoothly depending on
$\mu$, which is $O(\sqrt{\mu})$-shadowing the
chain $\sigma$:
$$
d(\gamma_\mu(t),C)\le \const\,\sqrt{\mu}.
$$
\item
Except for a small neighborhood $U$ of $M$ in $\calM$, $\gamma_\mu$ is  $O(\mu|\ln\mu|)$-shadowing $\sigma$:
\begin{equation}
\label{eq:d}
d(\gamma_\mu(t),C)\le \const\,\mu|\ln\mu|\quad\mbox{for}\quad \gamma_\mu(t)\in  \calM \setminus U.
\end{equation}
\item
The period of $\gamma_\mu$ is of order\footnote{The  notation means that the difference is   bounded  as $\mu\to 0$.}
\begin{equation}
\label{eq:period}
T_\mu\sim \sum_{i=0}^{n-1}\frac{|\ln\mu|}{\lambda(c_i)}.
\end{equation}
\end{itemize}
\end{thm}

\begin{rem}
The periodic orbit $\gamma_\mu$  has
$m$ pairs of multipliers  (eigenvalues of the linear Poincar\'e map) close  to the eigenvalues of $DF_n(c_0)$, and  $k-1$ pairs of hyperbolic multipliers $\rho,\rho^{-1}$ with  $|\rho|$ large of order $\mu^{-1}$. Thus $\gamma_\mu$
is always strongly unstable. If $DF_n(c_0)$ is hyperbolic, then $\gamma_\mu$ is a hyperbolic periodic orbit.
\end{rem}

The set $\cup_{0<\mu\le\mu_0}\gamma_\mu(\R)$
is a smooth invariant cylinder with
piece-wise smooth boundary $C\cup \gamma_{\mu_0}(\R)$.

If the chain $\sigma$ is negative, then  shadowing orbits exist
on $\Sigma_{\mu}$ with $\mu\in[-\mu_0,0)$.

A result similar to Theorem \ref{thm:period} holds for orbits shadowing nonperiodic heteroclinic chains.
Consider the skew product   of a finite subcollection $\calK$ of maps  $f\in \calF$.

\begin{thm}\label{thm:skew}
Let $\Lambda\subset \calK^\Z\times M$ be a compact hyperbolic invariant
set.  Take any orbit in
$\Lambda$ and let $\sigma=(\sigma_i)_{i\in\Z}$ be the corresponding
heteroclinic chain. Suppose that $\sigma$  is uniformly positive:
there is $\delta>0$ such that $a_i(\sigma)\ge \delta$ for all $i$.
There exists $\mu_0=\mu_0(\Lambda,\delta)$ such that for any $\mu\in
(0,\mu_0]$ there exists an orbit on $\Sigma_\mu$ which
$O(\sqrt{\mu})$-shadows the chain $\sigma$.
\end{thm}

When $M=\{z_0\}$ is a single hyperbolic equilibrium, a version of Theorem
\ref{thm:skew}  was proved in \cite{Bol-Mac:second} and used to study Poincar\'e second species solutions of the restricted circular 3 body problem. Then the scattering map is trivial, and so the nondegeneracy condition
for the heteroclinic chain does not appear. For $M=\{z_0\}$ and generic real eigenvalues,  an analog of Theorem
\ref{thm:skew} was announced in \cite{Tur-Shil2}.
The proof  appeared in \cite{Bol-Rab:revers}. In \cite{Tur-Shil2}
systems with discrete symmetries were studied. In \cite{Kaloshin},
regularity at $\mu=0$ of the cylinder formed by periodic orbits  was investigated
in relation to the problem of Arnold's diffusion.

In \cite{Bol-Rab:revers}
also global results on the existence of chaotic shadowing orbits were obtained by variational methods.
For a hyperbolic equilibrium with complex eigenvalues, shadowing via variational methods was done in \cite{Buf-Ser}.
We are not able to use global variational methods in the current setting.
although the proof of Theorem \ref{thm:period} has variational flavor.

The proof of Theorem \ref{thm:skew} is similar to that of Theorem \ref{thm:period}, but needs more work.
In order not to make the paper too long, we postpone this to a subsequent publication.
Also the existence of ``diffusion'' shadowing orbits with average speed along $M$ of order $|\ln\mu|^{-1}$ can be proved.
Note that this is much faster than in the problem of Arnold's diffusion, where (in the initially hyperbolic case)  the speed is of
order $O(\mu|\ln\mu|)$ \cite{Treschev:diff}. The reason is that we do not have the resonance gap problem.

Recently shadowing  chains of homoclinic orbits to a symplectic normally hyperbolic invariant manifold was studied in \cite{Delshams:shadow}  by the windows method. However, our situation is very different since the manifold $M$ is critical. In particular, in \cite{Delshams:shadow} the positivity condition does not appear.

As a corollary of Theorem \ref{thm:period}, we obtain a seemingly
more general bifurcation result. Consider a Hamiltonian
\begin{equation}
\label{eq:pert} H_\mu=H_0+\mu h+O(\mu^2)
\end{equation}
smoothly depending on the parameter $\mu$. Suppose $H_0$ satisfies the
conditions above, so it has a critical hyperbolic manifold $M\subset
\Sigma_0=H_0^{-1}(0)$ with real eigenvalues and (\ref{eq:Apm}) holds.
Let $\calF$ be the corresponding scattering map.

\begin{thm}\label{thm:pert}
Suppose $c_{i+1}=f_i(c_i)$ is a nondegenerate periodic  orbit of $\calF$ and
let $\sigma=(\sigma_i)$ be the corresponding periodic
heteroclinic chain of the flow $\phi^t_{H_0}$. Suppose that
$a_i(\sigma)h(c_i)<0$ for all $i$.
There exists $\mu_0>0$ such that for any $\mu\in (0,\mu_0]$ there
exists a periodic  orbit of the flow $\phi^t_{H_\mu}$ on $\Sigma_\mu=H_\mu^{-1}(0)$
which $O(\sqrt{\mu})$-shadows the  chain
$\sigma$. Moreover (\ref{eq:d})--(\ref{eq:period}) hold.
\end{thm}

A similar generalization of Theorem \ref{thm:skew} also holds.

If $h$ has constant sign on $\Sigma_0$, for example
$h|_{\Sigma_0}<0$, then Theorem \ref{thm:pert} immediately follows
from Theorem \ref{thm:period}. Indeed, in a compact subset of a
neighborhood of $\Sigma_0$  we can solve the equation $H_\mu(x)=0$
for
$$
\mu=\calH(x)=-\frac{H_0(x)}{h(x)}+\cdots
$$
and obtain  a Hamiltonian $\calH$ such that
$\calH^{-1}(\mu)=\Sigma_\mu$. Then the flows
$\phi^t_{H_\mu}|_{\Sigma_\mu}$ and $\phi^\tau_\calH|_{\Sigma_\mu}$
have the same trajectories, but with different time
parametrizations. Theorem \ref{thm:period} can be applied to the
flow $\phi^\tau_\calH$ which yields Theorem \ref{thm:pert}.

When $h$ changes sign,  one can define $\calH$ in
the domains $h>0$ and $h<0$, but not for $h=0$. Thus, in this case,
Theorem  \ref{thm:pert} does not follow from Theorem \ref{thm:period}.
However, the only place where there appear trajectories crossing the surface $h=0$  is in Corollary \ref{cor:conjugate}
 whose proof does not require
introduction of the Hamiltonian $\calH$. Thus the proof of Theorem \ref{thm:period} works for Theorem \ref{thm:pert}.
\qed

\medskip

The idea of the proof of Theorem \ref{thm:period} is variational.
We will construct a discrete action functional $\calA_\mu$, $\mu\in (0,\mu_0]$,
whose critical points correspond to  trajectories $\gamma_\mu$ on $\Sigma_\mu$ shadowing the heteroclinic chain $\sigma$.
The functional $\calA_\mu$  has a limit $\calA_0$ as $\mu\to0$
 and $\calA_\mu=\calA_0+O(\mu|\ln\mu|)$. A nondegenerate critical point of
the functional (\ref{eq:calA}) gives a nondegenerate critical point of $\calA_0$ and hence a nondegenerate critical point of $\calA_\mu$ for small $\mu$.

Construction of a functional $\calA_\mu$ continuous at $\mu=0$ is not evident, because $\gamma_\mu$ spends a long time of
order  $|\ln\mu|$ near $M$ and so, in some sense, the perturbation is singular at $\mu=0$.
The way out was found by Shilnikov \cite{Shilnikov} in the proof of the Shilnikov Lemma,
which is a version of the well known $\lambda$-lemma \cite{Katok}.
Shilnikov's method was used in \cite{Deng} to prove the strong $\lambda$-lemma.

The main result  of the present paper is Theorem \ref{thm:Shil_mu}
(generalization of the Shilnikov Lemma) which describes  solutions of a boundary value problem
for trajectories on $\Sigma_\mu$ near $M$.
It makes possible to construct a functional $\calA_\mu=\calA_0+O(\mu|\ln\mu|)$ and then prove Theorem \ref{thm:period}.
A weaker analog of Theorem \ref{thm:Shil_mu} was proved in \cite{Bol:DCDS}.

Theorem \ref{thm:Shil_mu}  was already used without proof in \cite{Bol-Neg:DCDS} to establish the existence of Poincar\'e second species solutions of the (nonrestricted) plane 3 body problem. So now the proof in \cite{Bol-Neg:DCDS}  is finally complete.
Application to the 3 body problem is briefly discussed in the next section.

 \section{Critical manifolds via Levi-Civita regularization in the 3 body problem}

\label{sec:Levi}

Consider the plane 3-body problem   with masses $m_1,m_2,m_3$. Suppose that $m_3$ is much larger than $m_1,m_2$:
$$
\frac{m_1}{m_3}=\mu\alpha_1,\quad \frac{m_2}{m_3}=\mu\alpha_2,\quad \alpha_1+\alpha_2=1,\quad \mu\ll 1.
$$
Let $q_1,q_2\in\R^2$ be positions of $m_1,m_2$ relative to $m_3$, and $p_1,p_2,p_3\in \R^2$ the
momenta. Setting $p_1+p_2+p_3=0$, we obtain the Hamiltonian
\begin{equation}
\label{eq:Hmu} H_\mu(q,p)=H_0(q,p)+\mu \left( \frac{|p_1+p_2|^2}2 - \frac{\alpha_1\alpha_2}{|q_1-q_2|}\right),
\end{equation}
where $q=(q_1,q_2)$, $p=(p_1,p_2)$. The unperturbed Hamiltonian
$$
H_0(q,p)=   \frac{|p_1|^2}{2 \alpha_1}+\frac{|p_2|^2}{2 \alpha_2} -\frac{ \alpha_1}{|q_1|}  -\frac{ \alpha_2}{|q_2|}.
$$
describes 2 uncoupled Kepler problems.

To regularize double collisions of $m_1,m_2$ at $\Delta=\{q_1=q_2\ne 0\}$, we identify $\R^2$ with $\C$ and perform the Levi-Civita symplectic transformation $g(x,y,\xi,\eta)=(q_1,q_2,p_1,p_2)$,
$$
q_1=x- \alpha_2\xi^2,\quad q_1=x+ \alpha_1\xi^2,\quad p_1= \alpha_1 y-\frac{\eta}{2\bar\xi},\quad p_2= \alpha_2y+\frac{\eta}{2\bar\xi}.
$$
The map $g$ is  a double covering undefined at $\xi=0$ which corresponds to double collisions at $\Delta$.
We fix energy $E$ and set
\begin{eqnarray}
\label{eq:calH}
&\calH_\mu^E(x,y,\xi,\eta)=|\xi|^2(H_\mu\circ g-E)\\
&=\frac{|\eta|^2}{8 \alpha_1\alpha_2}-|\xi|^2\left( E+\frac{ \alpha_1}{| \alpha_2\xi^2-x|} +\frac{ \alpha_2}{|
\alpha_1\xi^2+x|}-\frac{(1+\mu)|y|^2}2 \right)+\mu\alpha_1\alpha_2.\nonumber
\end{eqnarray}
Denote  $\Sigma_\mu^E=H_\mu^{-1}(E)$ and $\Gamma_\mu^E=(\calH_\mu^E)^{-1}(0)$. Since $g(\Gamma_\mu^E)=\Sigma_\mu^E$, the
map $g$ takes orbits of  the flow $\phi^\tau_{\calH_\mu^E}$  on $\Gamma_\mu^E$ to orbits of  the flow $\phi^t_{H_\mu}$ on $\Sigma_\mu^E$. The time parametrization is changed: the new time is given by $d\tau=|\xi|^2\,dt$.

The singularity at $\Delta$ disappeared: the regularized Hamiltonian $\calH_\mu^E$ is smooth on
$$
 \calM =\{(x,y,\xi,\eta):x\ne \alpha_2\xi^2,\;x\ne-\alpha_1\xi^2\}
$$
which means excluding collisions of $m_1$ and $m_2$ with $m_3$. Double collisions of $m_1$ and $m_2$ correspond to
$\xi=\eta=0$. For $\mu=0$, the Hamiltonian
$$
\label{eq:calH_0} \calH_0^E(x,y,\xi,\eta)=\frac{|\eta|^2}{8
\alpha_1\alpha_2}-|\xi|^2\Big(E+\frac{1}{|x|}-\frac{|y|^2}{2}\Big)+O(|\xi|^4)
$$
has a   normally hyperbolic symplectic critical manifold
$$
 M_E= \{(x,y,0,0):  \frac12|y|^2-\frac{1}{|x|}<E\}
$$
with real semisimple eigenvalues
$$
\pm\sqrt{\frac{1}{2\alpha_1\alpha_2}\bigg(E+\frac{1}{|x|}-\frac{|y|^2}{2}\bigg)}.
$$

For $\mu=0$, collision orbits of $m_1,m_2$ (pairs of arcs of Kepler orbits starting and ending at $\Delta$) with energy $E$ correspond to
trajectories of $\phi_{\calH_0^E}^\tau$ asymptotic to $M_E$, and chains of collision orbits with continuous total momentum $y=p_1+p_2$ correspond to chains of heteroclinic orbits. For small $\mu>0$, orbits of the 3 body problem  with energy $E$ passing $O(\mu)$-close to the singular set $\Delta$  correspond to orbits of the flow $\phi^\tau_{\calH_\mu^E}$  on the level $\Gamma_\mu^E$ passing $O(\sqrt{\mu})$-close to  $ M_E$.

The Hamiltonian (\ref{eq:calH}) has the   form (\ref{eq:pert}):
$$
\calH_\mu^E=\calH_0^E+\mu h,
$$
where $h|_{M_E} =\alpha_1\alpha_2>0$. Thus we are in the situation of Theorem \ref{thm:pert}. In \cite{Bol-Neg:DCDS} many
nondegenerate periodic collision chains to $M_E$ were obtained. Then for small $\mu>0$ Theorem \ref{thm:pert}
implies the existence of many periodic  almost collision solutions of the 3 body problem. Such solutions were named by Poincar\'e second species solutions. See \cite{Bol-Neg:DCDS} for details.

The plan of the paper is as follows. In section 3 we represent the stable and unstable manifolds by generating functions. In section 4  different versions of local connection theorems are formulated.
The proofs are given in section 5.
In section 6 relations between the generating functions of the scattering map
 and of the stable and unstable manifolds are discussed. In section 7 trajectories shadowing heteroclinic chains are
represented by critical points of a discrete action  functional, and then
Theorem \ref{thm:period} is proved.

\section{Generating functions of the stable and unstable manifolds}

In this section it does not matter if the eigenvalues
of critical points in $M$ are real or complex: we only need the critical manifold
$M$ to be symplectic and normally hyperbolic.

Take an open
set $V\Subset M$  with symplectic
coordinates $z=(x,y)\in\R^{2m}$ and
identify $V$ with a domain in $\R^{2m}$.
If $V$ is small enough, the stable and unstable bundles $E^\pm|_V$ are trivial over $V$.
Hence a tubular neighborhood $U$ of $V$ in $ \calM$ can be
identified with
$$
U\cong V\times B_r\times B_r =\{(z,q,p):z\in V,\; q,p\in B_r\},\qquad
B_r=\{q\in \R^{k}:|q|\le r\},
$$
in such a way  that $V\cong V\times(0,0)$ and for $z\in M$,
$$
E_z\cong\R^k\times \R^k,\quad E_z^+\cong\R^{k}\times\{0\}, \quad E_z^-\cong\{0\}\times\R^{k}.
$$
By the generalized Darboux Theorem (see \cite{Sympl}), we can assume that the coordinates  in $U$ are symplectic:
$$
\omega|_U=dy\wedge dx+dp\wedge dq.
$$
Then for $\xi=(\xi_+,\xi_-)$ and $\eta=(\eta_+,\eta_-)$ in
$E_z$,
\begin{equation}\label{eq:omega1}
\omega(\xi,\eta)=\langle \xi_-,\eta_+\rangle- \langle \eta_-,\xi_+\rangle.
\end{equation}

Since the local stable and unstable manifolds $W_\loc^\pm(V)$ are tangent to $E^\pm|_V$, they are graphs
\begin{equation}
\begin{array}{l}
\label{eq:W+-} W^+_\loc(V)=\{(z,q,p):z\in V,\;q\in B_r,\; p=f_+(z,q)\},\\
W^-_\loc(V)=\{(z,q,p):z\in V,\;p\in B_r,\; q=f_-(z,p)\},
\end{array}
\end{equation}
where
$$
f_+(z,q)=O_2(q),\quad  f_-(z,p)=O_2(p).
$$

\begin{rem}
$O_2(q)$ means a function of the form $\sum_{|i|=2}a_{i}(z,q)q^i$ with smooth
coefficients. For $i\in\Z_+^{k}$ we write $|i|=i_1+\cdots+i_{k}$.
\end{rem}

Take a smaller open set $V_0\Subset V$. For any $z_0\in V_0$ the local stable and unstable manifolds  are given by
$W^\pm_\loc(z_0)= \psi_\pm(z_0,B_r)$, where
\begin{eqnarray*}
\psi_+(z_0,q)=(g_+(z_0,q),q,h_+(z_0,q))=(z_0,q,0)+O_2(q),\\
 \psi_-(z_0,p)=(g_-(z_0,p),h_-(z_0,p),p)=(z_0,0,p)+O_2(p).
\end{eqnarray*}
and
\begin{eqnarray*}
h_+(z_0,q)=f_+(g_+(z_0,q),q),\quad h_-(z_0,p)=
f_-(g_-(z_0,p),p).
\end{eqnarray*}
For $z_0\in V_0$ and $q_+,p_-\in B_r$ let
\begin{equation}
\begin{array}{l}
\label{eq:gamma+-}
\gamma_+:[0,+\infty)\to W^+_\loc(z_0),\quad
\gamma_+(t)=\phi^t\circ \psi_+(z_0,q_+),\\
 \gamma_-:(-\infty,0]\to W_\loc^-(z_0),\quad
\gamma_-(t)=\phi^t\circ \psi_-(z_0,p_-),
\end{array}
\end{equation}
be the trajectories asymptotic to $z_0$ as $t\to \pm\infty$. Then
\begin{eqnarray*}
\gamma_+(0)=\psi_+(z_0,q_+)=(z_+,q_+,p_+),\\
\gamma_-(0)=\psi_-(z_0,p_-)=(z_-,q_-,p_-).
\end{eqnarray*}

We will represent $W_\loc^\pm(z_0)$ by generating functions
as follows.

\begin{prop}\label{prop:gen}
There exist smooth functions
\begin{equation}
\begin{array}{l}
\label{eq:Spm} S_+(x_+,y_0,q_+)=\langle x_+,y_0\rangle +O_2(q_+),\\
S_-(x_0,y_-,p_-)=\langle x_0,y_-\rangle +O_2(p_-),
\end{array}
\end{equation}
on open sets in $\R^m\times\R^m\times \R^{k}$ such that for any $z_0=(x_0,y_0)\in V_0$ and
$A_\pm=(x_\pm,y_\pm,q_\pm,p_\pm)\in U$,
\begin{eqnarray}
A_+\in W_\loc^+(z_0)\;\Leftrightarrow\; p_+=\frac{\partial
S_+}{\partial q_+},\quad y_+=\frac{\partial S_+}{\partial x_+},
\quad x_0= \frac{\partial S_+}{\partial y_0},\label{eq:pi+}\\
A_-\in W_\loc^-(z_0)\;\Leftrightarrow\; q_-=\frac{\partial
S_-}{\partial p_-},\quad x_-=\frac{\partial S_-}{\partial y_-},
\quad y_0= \frac{\partial S_-}{\partial x_0}.\label{eq:pi-}
\end{eqnarray}
\end{prop}

Equivalently,
\begin{eqnarray}\label{eq:deltaS+}
 d S_+(x_+,y_0,q_+)&=&p_+\, d q_++y_+\, d x_++x_0\, d y_0,\\
 d S_-(x_0,y_-,p_-)&=&q_-\, d p_-+x_-\, d
y_-+y_0\, d x_0. \label{eq:deltaS-}
\end{eqnarray}
In particular,
$$
(x_+,q_+)\to S_+(x_+,y_0,q_+),\quad (y_-,p_-)\to S_-(x_0,y_-,p_-)
$$
are the generating functions of the Lagrangian manifolds
$W^+_\loc(y=y_0)$ and $W^-_\loc(x=x_0)$.

\proof
Let
$$
J_+(z_0,q_+)=\int_{\gamma_+}\alpha,\quad J_-(z_0,p_-)=\int_{\gamma_-}\alpha,\qquad \alpha=y\,dx+p\,dq.
$$
be the Maupertuis actions of the asymptotic trajectories $\gamma_\pm$.
The first variation formula \cite{Arnold} gives
\begin{eqnarray}\label{eq:var+}
d J_+(z_0,q_+)&=&y_0\, d x_0-y_+\, d x_+ - p_+\, d q_+,\\
d J_-(z_0,p_-)&=& y_-\, d x_- +
p_-\, d q_- - y_0\, d x_0. \label{eq:var-}
\end{eqnarray}

Equations (\ref{eq:var+})--(\ref{eq:var-})    imply that
\begin{eqnarray*}
z_0\to z_+=g_+(z_0,q_+)=z_0+O_2(q_+),\\
z_0\to z_-=g_-(z_0,p_-)=z_0+O_2(p_-),
\end{eqnarray*}
are symplectic  maps which are close to identity. We represent them
by appropriate generating functions \cite{Arnold}. Let
\begin{eqnarray}
g_+(z_0,q_+)=(X_+(z_0,q_+),Y_+(z_0,q_+)),\label{eq:g+-}\\
g_-(z_0,p_-)=(X_-(z_0,p_-),Y_-(z_0,p_-)).
\end{eqnarray}
Set
\begin{eqnarray*}
S_+(x_+,y_0,q_+)=\langle y_0,x_0\rangle -J_+(z_0,q_+),
\end{eqnarray*}
where $ x_0(x_+,y_0,q_+)$ is a  solution of the  equation
\begin{equation}
\label{eq:X+}
x_+=X_+(x_0,y_0,q_+)=x_0+O_2(q_+).
\end{equation}
Similarly, set
$$
S_-(x_0,y_-,p_-)=\langle y_-,x_-\rangle+\langle p_-,q_-\rangle
-J_-(z_0,p_-),
$$
where $(z_-,q_-,p_-)=\psi_-(z_0,p_-)$ and $y_0(x_0,y_-,p_-)$ is a solution of  the equation
\begin{equation}
\label{eq:Y-}
y_-=Y_-(x_0,y_0,p_-)=y_0+O_2(p_-).
\end{equation}

By (\ref{eq:var+})--(\ref{eq:var-}), the functions $S_\pm$ satisfy
(\ref{eq:deltaS+})--(\ref{eq:deltaS-}). \qed

\medskip

Next we combine asymptotic orbits $\gamma_\pm$ in one  curve $\gamma_+\cdot \gamma_-$
with reflection from $M$ at  $z_0$.  If $r>0$ is small enough,
for any\footnote{The notation $(x_+, y_-)\in V_0$ makes sense because
we identified $V_0$ with a domain in $\R^{2m}$.} $(x_+, y_-)\in V_0$ and $q_+,p_-\in B_r$
we can solve equations (\ref{eq:X+})--(\ref{eq:Y-})  for
\begin{equation}
\label{eq:zeta}
z_0=\zeta(Z)=(x_+,y_-)+O_2(q_+,p_-),\qquad Z=(x_+,y_-,q_+,p_-).
\end{equation}

\begin{prop}\label{prop:reflect}
Suppose $r>0$ is sufficiently small.
Then for any  $Z=(x_+, y_-,q_+,p_-)\in V_0\times B_r\times B_r$:
\begin{itemize}
\item There exist $z_0\in V$, $x_-,y_+\in\R^m$, and $q_-,p_+\in \R^{k}$ such that
\begin{equation}\label{eq:reflect}
\begin{array}{l}
A_+=(x_+,y_+,q_+,p_+)=\psi_+(z_0,q_+)\in W_\loc^+(z_0),\\
A_-=(x_-,y_-,q_-,p_-)=\psi_-(z_0,p_-)\in W_\loc^-(z_0).
\end{array}
\end{equation}
\item
The relation $A_+\to A_-$ is symplectic: there is a smooth generating function
\begin{equation}
\label{eq:L2} L(Z)=\langle x_+,y_-\rangle +O_2(q_+,p_-).
\end{equation}
such that (\ref{eq:reflect}) is equivalent to
    \begin{equation}
\label{eq:L}
 d L(Z)= y_+\, d x_+ + x_-\, d y_- + p_+\, d q_+ +q_-\, d p_-.
\end{equation}
\end{itemize}
\end{prop}

\proof
Consider the function
$$
F(z_0,Z)=S_+(x_+,y_0,q_+)+
S_-(x_0,y_-,p_-)-\langle x_0,y_0\rangle.
$$
Then (\ref{eq:pi+})--(\ref{eq:pi-}) imply that
$A_+\in W_\loc^+(z_0)$ and $A_-\in W_\loc^-(z_0)$ iff
\begin{equation}
\label{eq:F}
\frac{\partial F}{\partial z_0}=0,\quad y_+=\frac{\partial F}{\partial x_+},\quad
x_-=\frac{\partial F}{\partial y_-},\quad p_+=\frac{\partial
F}{\partial q_+}, \quad q_-=\frac{\partial F}{\partial p_-}.
\end{equation}
We have
$$
z_0=\zeta(Z)\;\Leftrightarrow\;  \frac{\partial F}{\partial z_0}=0.
$$
Define the generating function $L$ by
\begin{equation}
\label{eq:Crit}
L(Z)=F(\zeta(Z),Z)=\Crit_{z_0}F(z_0,Z)
\end{equation}
which means taking the nondegenerate critical value with respect to $z_0$.
Then (\ref{eq:F}) implies (\ref{eq:L}).
 \qed

\begin{rem}
The generating function $L$ does not satisfy the twist condition. Indeed, a computation gives
$$
\left(
\begin{array}{cc}
\frac{\partial^2 L}{\partial x_+\partial y_-} &\frac{\partial^2 L}{\partial x_+\partial p_-}\\ \frac{\partial^2 L}{\partial
q_+\partial y_-} & \frac{\partial^2 L}{\partial q_+\partial p_-}
\end{array}
\right)
=
\left(
\begin{array}{c}
\frac{\partial x_0}{\partial x_+}\\ \frac{\partial x_0}{\partial q_+}
\end{array}
\right)
\left(
\frac{\partial y_0}{\partial y_-}, \frac{\partial y_0}{\partial p_-}
\right)
$$
Hence the rank of this matrix is $m$. Equations (\ref{eq:L}) do not define a
map $ A_+\to A_-$. The correspondence $A_+\to A_-$ is a symplectic relation,
i.e.\ a Lagrangian submanifold in $ \calM\times  \calM $.
\end{rem}

\section{Local connection}

In this section we formulate several connection theorems
describing the behavior of trajectories of the Hamiltonian system  near
the critical manifold $M$. In the rest of the paper we  assume
(\ref{eq:Apm}). In the generic case (\ref{eq:generic}) the results
are similar, but they  will be published elsewhere.

By (\ref{eq:d2H2}), in the coordinates $(z,q,p)$ in a tubular neighborhood $U\cong V\times B_r\times
B_r$ of $V\Subset M$, the Hamiltonian has the form
\begin{equation}
\label{eq:H|U} H|_U=H(z,q,p)=-\lambda(z)\langle p,q\rangle +O_3(p,q) .
\end{equation}
The corresponding Hamiltonian system is
\begin{eqnarray*}
\dot z&=&O_2(p,q),\\ \dot q&=&\frac{\partial H}{\partial p}=-\lambda(z)q+O_2(p,q),\\
 \dot p&=&-\frac{\partial H}{\partial q}=\lambda(z)p+O_2(p,q).
\end{eqnarray*}
The limit directions (\ref{eq:dir}) of the asymptotic orbits
(\ref{eq:gamma+-}) are
\begin{eqnarray*}
&\vbf_+(\gamma_+)=(0,v_+,0),\qquad &v_+(z_0,q_+)=\lim_{t\to+\infty}e^{\lambda(z_0)t}q(t)=q_++O_2(q_+),\\
&\vbf_-(\gamma_-)=(0,0,v_-),\qquad &v_-(z_0,p_-)=\lim_{t\to-\infty}e^{-\lambda(z_0)t}p(t)=p_-+O_2(p_-).
\end{eqnarray*}
By (\ref{eq:omega1}), the symplectic angle of the concatenation $\gamma_+\cdot\gamma_-$ at $z_0$ is
\begin{equation}
\label{eq:omega}
\omega(\vbf_+(\gamma_+),\vbf_-(\gamma_-))=- \langle v_+(z_0,q_+),v_-(z_0,p_-)\rangle =-\langle q_+,p_-\rangle+O_3(q_+,p_-).
\end{equation}

There are two main versions of connection theorems: for fixed time and for fixed energy.

\begin{thm}[Fixed time connection]
\label{thm:Shil}
Suppose that $r>0$ is small enough.
For any  $Y=(z_0,q_+,p_-)\in V_0\times  B_r\times B_r$ and $T\ge 1$:
\begin{itemize}
\item There exists a unique solution
\begin{equation}
\label{eq:gamma}
\gamma(t)=(z(t),
 q(t),p(t))\in V\times B_{r}\times B_{r},\qquad t \in [- T, T],
\end{equation}
satisfying the initial--boundary conditions
\begin{equation}
\label{eq:BC0} z(0) = z_0,\quad  p(T)=p_-,\quad q(-T)=q_+.
\end{equation}
\item
$\gamma$ smoothly depends on  $(Y,T)\in V_0\times B_r\times
B_r\times [1,+\infty)$.
\item
$\gamma(t)$ converges to
$\gamma_+(t+T)$   on $[-T,0]$ and
to $\gamma_-(t-T)$   on $[0,T]$ as $T\to\infty$:
\begin{eqnarray}
\label{eq:bolP} \gamma(t)=\gamma_+(t+T)+\gamma_-(t-T)-(z_0,0,0)
 +e^{-\lambda(z_0)T}O(r^2).
\end{eqnarray}
Thus $\gamma([-T,T])$ converges to the concatenation $\gamma_+\cdot\gamma_-$.
\item  Let
\begin{eqnarray}
\label{eq:a+-}
\gamma(\mp T)=A_\pm=(z_\pm,q_\pm,p_\pm),\quad \gamma(0)=(z_0,q_0,p_0).
\end{eqnarray}
Then
\begin{eqnarray}
z_+&=&g_+(z_0,q_+)+Te^{-2\lambda(z_0)T}O(r^2),\label{eq:z+}\\
z_-&=&g_-(z_0,p_-)+Te^{-2\lambda(z_0)T}O(r^2),\label{eq:z-}\\
p_+&=& h_+(z_0,q_+)+e^{-2\lambda(z_0)T}(p_-+O(r^2)),\\
q_-&=& h_-(z_0,p_-)+e^{-2\lambda(z_0)T}(q_++O(r^2)),\label{eq:q-}\\
q_0&=&e^{-\lambda(z_0)T}
v_+(z_0,q_+)+e^{-2\lambda(z_0)T}O(r^2),\label{eq:q(0)}\\
p_0&=&e^{-\lambda(z_0)T}v_-(z_0,p_-)+e^{-2\lambda(z_0)T}O(r^2).\label{eq:p(0)}
\end{eqnarray}
\end{itemize}
\end{thm}

\begin{rem}
Here $O(r^2)$ means a function $f(z_0,q_+,p_-,T)$ on $V\times
B_r\times B_r$, depending also on $T\ge 1$, such that
\begin{equation}
\label{eq:norm}
\|f\|_{C^1(V\times B_r\times
B_r)}=\sup_{V\times B_r\times
B_r}\max\bigg\{|f|,\bigg|\frac{\partial f}{\partial z_0}\bigg|,r \bigg|\frac{\partial f}{\partial q_+}\bigg|,r\bigg|\frac{\partial f}{\partial p_-}\bigg|\bigg\} \le Cr^2
\end{equation}
with   $C$ independent of $r$ and $T$. Thus the norms of the
derivatives with respect to $q_+,p_-\in B_r$ are taken with weight
$r$. Equivalently, (\ref{eq:norm}) is the $C^1$ norm of the function $f(z_0,r\hat q_+,r\hat p_-,T)$
on $V\times B_1\times B_1$.
\end{rem}

When $M=\{z_0\}$ is a single equilibrium with equal eigenvalues,
Theorem \ref{thm:Shil} was proved in \cite{Bol-Mac:second}.  When $M$ is a single equilibrium
with generic real eigenvalues, an analog of Theorem \ref{thm:Shil}  can be deduced from the strong $\lambda$-lemma \cite{Deng},
see \cite{Bol-Rab:revers}.

With minor modifications  Theorem \ref{thm:Shil}
holds also for non-Hamiltonian and non-autonomous systems.
Now we will use the Hamiltonian structure.
For large $T$, we solve  (\ref{eq:z+})--(\ref{eq:z-}) for
\begin{equation}
z_0=\zeta_T(Z)=\zeta(Z)+ Te^{-2\lambda(\zeta)T}O(r^2),\qquad Z=(x_+,y_-,q_+,p_-),\label{eq:last}
\end{equation}
where $\zeta$ is the function (\ref{eq:zeta}). We obtain

\begin{thm} \label{thm:x+y-}
Let $T_0>0$ be sufficiently large. For every $T\ge T_0$
and $Z=(x_+,y_-,q_+,p_-)\in V_0\times  B_r\times B_r$:
\begin{itemize}
\item
There exists a solution (\ref{eq:gamma}) satisfying the boundary conditions
\begin{equation}
\label{eq:bound0}
x(-T)=x_+,\quad y(T)=y_-,\quad q(-T)=q_+,\quad p(T)=p_-.
\end{equation}
\item
The relation $A_+\to A_-$ between the points (\ref{eq:a+-})
is symplectic: there exists
a smooth generating function $L_T(Z)$ such that
$$
d L_T(Z)= y_+\, d x_+ + x_-\, d y_- +
p_+\, d q_+ + q_-\, d p_-.
$$
\item As $T\to+\infty$, the generating function has the asymptotics
\begin{equation}\label{eq:ST}
L_T(Z)=L(Z)
+e^{-2\lambda(\zeta)T}(\langle
q_+,p_-\rangle+TO(r^3))
\end{equation}
where  $L$ is the generating function (\ref{eq:L}).
\end{itemize}
\end{thm}

Since generating functions are defined up to a constant, the equality (\ref{eq:ST}) is modulo a constant.
The symplectic relation $A_+\to A_-$ has  a smooth limit as $T\to
+\infty$. This is true because of  a right choice of the boundary
conditions which is motivated by the Shilnikov Lemma \cite{Shilnikov}.
For small $r$, the generating function $L_T$ satisfies the twist condition, but the twist   is
exponentially small for $T\to+\infty$.

\medskip

Next we formulate the fixed energy version of the connection theorem.
Fix arbitrary $\nu,\kappa\in (0,1)$ and let
\begin{equation}
\label{eq:DQ}
D_r=B_{r}\setminus B_{\nu r},\quad
Q_r=\{(q_+,p_-)\in D_r\times D_r: \langle q_+,p_-\rangle\le -\kappa
r^2\}.
\end{equation}

\begin{thm}[Fixed energy connection]
\label{thm:Shil_mu}
Let $r>0$ and $\mu_0>0$ be sufficiently small.
Then for any $\mu\in (0,\mu_0]$  and $Y=(z_0,q_+,p_-)\in V_0\times Q_r$:
\begin{itemize}
\item There exist
\begin{equation}
\label{eq:T}
T=\frac{|\ln\mu|+\ln(-\lambda(z_0)\langle
v_+(z_0,q_+),v_-(z_0,p_-)\rangle)}{2\lambda(z_0)}
+O(\sqrt{\mu})
\end{equation}
and a unique solution (\ref{eq:gamma}) on $\Sigma_\mu=H^{-1}(\mu)$ satisfying (\ref{eq:BC0}).
\item $\gamma$ and $T$ smoothly depend on  $(Y,\mu)\in V_0\times Q_r\times
(0,\mu_0]$.
\item $\gamma$ converges to the concatenation $\gamma_+\cdot\gamma_-$  as $\mu\to 0$:
\begin{eqnarray}
\label{eq:concat} \gamma(t)=\gamma_+(t+T)+ \gamma_-(t-T)-(z_0,0,0)+O(\sqrt{\mu}).
\end{eqnarray}
\item  The points (\ref{eq:a+-}) satisfy
\begin{eqnarray}
z_+&=&g_+(z_0,q_+)+O(\mu|\ln\mu|),\label{eq:z_+}\\
z_-&=&g_-(z_0,p_-)+O(\mu|\ln\mu|),\label{eq:z_-}\\
p_+&=&h_+(z_0,q_+)- \frac{ \mu p_-}{\lambda(z_0)\langle
q_+,p_-\rangle}
+ O(\mu),\label{eq:p_+}\\
q_-&=& h_-(z_0,p_-)-  \frac{\mu q_+}{\lambda(z_0)\langle
q_+,p_-\rangle}+  O(\mu),\label{eq:p_-}\\
q_0&=&\sqrt{\mu}\,v_+(z_0,q_+)+ O(\mu),\\
p_0&=&\sqrt{\mu}\,v_-(z_0,p_-)+O(\mu ).
\end{eqnarray}
\end{itemize}
\end{thm}

For the case when $M$ is a single equilibrium, Theorem \ref{thm:Shil_mu}
was obtained in \cite{Bol-Mac:second}. A version of Theorem \ref{thm:Shil_mu} was used without proof in
\cite{Bol-Neg:DCDS}.

\begin{rem}
Here $O(\mu)$ or $O(\mu|\ln\mu|)$ means a function $f$ on $V\times
Q_r$, depending also on $\mu\in (0,\mu_0]$, such that
$$
\|f\|_{C^1(V\times
Q_r)}  \le C\mu\quad \mbox{or}\quad \|f\|_{C^1(V\times
Q_r)}  \le C\mu|\ln\mu|,
$$
where the constant is independent of $r$ and $\mu$. The  $C^1$ norm is weighted as in (\ref{eq:norm}).
Hence the second terms in (\ref{eq:p_+})--(\ref{eq:p_-}) are not $O(\mu)$.
They provide nontrivial twist in the Poincar\'e map,
see Remark \ref{rem:twist}.
\end{rem}

\begin{rem}
By (\ref{eq:omega}), for small $r$, $(q_+,p_-)\in Q_r$ implies
$\omega(\vbf_+(\gamma_+),\vbf_-(\gamma_-))>0$.
Thus the concatenation of $\gamma_+$ and $\gamma_-$ at $z_0$ is positive (see Definition \ref{def:pos}).
This explains how the positivity condition appears in Theorem \ref{thm:period}.
If we replace the set $Q_r$ by
$$
\{(q_+,p_-)\in D_r\times D_r: \langle q_+,p_-\rangle\ge \kappa
r^2\},
$$
then the concatenation of $\gamma_+$ and $\gamma_-$ at $z_0$ is negative, and
the connecting solution $\gamma$ exists for $\mu\in[-\mu_0,0)$.
\end{rem}

\begin{rem}
For simplicity we fixed $\kappa>0$ in (\ref{eq:DQ}). In fact Theorem
\ref{thm:Shil_mu} can be improved to include $\kappa=C\mu^{1/3}$
with $C>0$ sufficiently large constant. However, we do not need this
for our purposes.
\end{rem}

Let us deduce Theorem \ref{thm:Shil_mu} from Theorem \ref{thm:Shil}.
By (\ref{eq:q(0)})--(\ref{eq:p(0)})
and (\ref{eq:H|U}), on the connecting trajectory $\gamma$ in Theorem
\ref{thm:Shil},
$$
\label{eq:H} H|_\gamma=H(\gamma(0))=-\lambda(z_0)e^{-2\lambda(z_0)T}\langle
v_+(z_0,q_+),v_-(z_0,p_-)\rangle + e^{-3\lambda(z_0)T}O(r^3).
$$
To find  $\gamma$ on $\Sigma_\mu$, we solve the
equation $H|_\gamma=\mu$  for $T$. For $(q_+,p_-)\in Q_r$ and small
$\mu>0$ we obtain
$$
e^{-2\lambda(z_0)T}=-\frac{\mu+O(\mu^{3/2})}{\lambda(z_0)\langle
v_+(z_0,q_+),v_-(z_0,p_-)\rangle}>0.
$$
This implies (\ref{eq:T}) and Theorem \ref{thm:Shil_mu} follows easily.
In the next section we give an independent proof of Theorem \ref{thm:Shil_mu}.

\medskip

Solving   (\ref{eq:z_+})--(\ref{eq:z_-}) for $z_0$,
we obtain a symplectic version of the fixed energy connection theorem.

\begin{thm}
\label{thm:xy_mu}
Let $r>0$ and $\mu_0>0$  be sufficiently small.
Then for any $\mu\in (0,\mu_0]$ and $Z=(x_+,y_-,q_+,p_-)\in V_0\times Q_r$:
\begin{itemize}
\item
There exist
$$
z_0=\zeta_\mu(Z)=\zeta(Z)+O(\mu|\ln\mu|)
$$
and a solution (\ref{eq:gamma}) on $\Sigma_\mu$ satisfying boundary conditions (\ref{eq:bound0})
with $T=T_\mu(Z)$ as in (\ref{eq:T}).
\item 
$\gamma$ and $T$ smoothly depend on  $(Z,\mu)\in V_0\times Q_r\times
(0,\mu_0]$.
\item
The relation $A_+\to A_-$ between the points (\ref{eq:a+-}) is given by
\begin{eqnarray}
y_+&=&\frac{\partial L}{\partial x_+}+O(\mu|\ln\mu|),\\
x_-&=&\frac{\partial L}{\partial y_-}+O(\mu|\ln\mu|),\\
p_+&=&\frac{\partial L}{\partial q_+}- \frac{\mu p_-}{\lambda(\zeta)\langle
q_+,p_-\rangle}
+ O(\mu),\\
q_-&=& \frac{\partial L}{\partial p_-}-  \frac{\mu q_+}{\lambda(\zeta)\langle
q_+,p_-\rangle}+  O(\mu).
\end{eqnarray}
\item
The generating function of the symplectic relation $A_+\to A_-$ has the form
\begin{equation}
\label{eq:S_mu}
R_\mu(Z)=L(Z)-\frac{\mu\ln|\langle q_+,p_-\rangle|}{\lambda(x_-,y_+)}+O(\mu|\ln\mu|).\end{equation}
\end{itemize}
\end{thm}

Here $O(\mu|\ln\mu|)$ means a function with
$\|f\|_{C^2(V\times Q_r)}\le C\mu|\ln\mu|$,
where  $C$ is independent of $r,\mu$, and the norm is weighted as in (\ref{eq:norm}).

Theorem \ref{thm:xy_mu} follows  from Theorem \ref{thm:Shil_mu} and the implicit function theorem.
Conversely, Theorem \ref{thm:Shil_mu} can be deduced from Theorem \ref{thm:xy_mu}.
 We prove Theorems \ref{thm:x+y-} and \ref{thm:xy_mu} in the next section.
The proof of Theorem \ref{thm:Shil} is similar and we skip it.

The relation  $A_+\to A_-$ is restricted to the contact manifold
$\Sigma_\mu$. To get a symplectic map, we   take symplectic cross sections
\begin{equation}
\label{eq:Nmu}
\begin{array}{l}
N_\mu^+=\{(z_+,q_+,p_+)\in U\cap\Sigma_\mu:q_+\in S_r\},\\
N_\mu^-=\{(z_-,q_-,p_-)\in U\cap\Sigma_\mu:p_-\in S_r\},
\end{array}
\end{equation}
where $S_r=\partial B_r$ is a sphere.

\begin{cor}
The restriction of the function $R_\mu$ to the set
$$
E_r=\{Z=(x_+,y_-,q_+,p_-)\in V\times Q_r: q_+,p_-\in S_r\}
$$
is the generating function of the local Poincar\'e  map
$P_\mu:N_\mu^+\cap O^+\to N_\mu^-\cap O^-$:
$$
d R_\mu(Z)= y_+\, d x_+ + x_-\, d y_- + p_+\, d
q_+ + q_-\, d p_-.
$$
\end{cor}

Here $O^\pm$ are open sets in $U$. We introduce local symplectic coordinates $x_\pm,y_\pm,\xi_\pm,\eta_\pm$
on $N_\mu^\pm\cap O^\pm$  such that
\begin{equation}
\label{eq:local-sympl}
\begin{array}{l}
(y_+\,dx_++p_+\,dq_+)|_{N_\mu^+}=y_+\,dx_++\eta_+\,d\xi_+,\\
(x_-\,dy_-+q_+\,dp_-)|_{N_\mu^-}=x_-\,dy_-+\xi_-\,d\eta_-.
\end{array}
\end{equation}
Then
$$
\calR_\mu(x_+,y_-,\xi_+,\eta_-)=R_\mu(x_+,y_-,q_+(\xi_+),p_-(\xi_-))
$$
is the generating function of the coordinate representation of the
Poincar\'e  map $P_\mu:(x_+,y_+,\xi_+,\eta_+)\to
(x_-,y_-,\xi_-,\eta_-)$:
 $$
\calR_\mu(x_+,y_-,\xi_+,\eta_-)= y_+\, d x_+ + x_-\, d y_- +
\eta_+\, d \xi_+ + \xi_-\, d \eta_-.
$$

The coordinates $x_\pm,y_\pm,\xi_\pm,\eta_\pm$ on $N_\mu^\pm$  are defined as follows.
Choose local coordinates on the sphere $S_r$, for example given by a stereographic projection. Then
$q_+=q_+(\xi_+)\in S_r$ and $p_-=p_-(\eta_-)\in S_r$, where
$\xi_+,\eta_-\in\R^{k-1}$.  Set
\begin{equation}
\label{eq:local}  \eta_+=p_+\cdot Dq_+(\xi_+),\quad \xi_-=q_-\cdot
Dp_-(\eta_-).
\end{equation}
Then $(z_\pm,\xi_\pm,\eta_\pm)$ determine $(z_\pm,q_\pm,p_\pm)$ and so they are local coordinates on $N_\mu^\pm$.

Indeed, let $(z_+,q_+(\xi_+),p_+)\in N_\mu^+$.  The orthogonal
projection $\bar p_+\perp q_+$ of $p_+$ to $T_{q_+}S_r$ is
determined by $\eta_+=\bar p_+\cdot Dq_+(\xi_+)$. Then $p_+=c q_+ +\bar p_+$, where the scalar
$c$  is the solution of the equation
$$
H(z_+,q_+,p_+) =\lambda(z_+)\langle q_+,p_+\rangle+O(r^3)=
\lambda(z_+)r^2c+O(r^3)=\mu.
$$

\begin{rem}\label{rem:twist}
The generating function $L|_{E_r}$ does not satisfy the twist condition, but the
function $R_\mu|_{E_r}$ does, with the twist in $q_+,p_-$ of order $\mu$.
Thus we  are in the situation of the so called anti-integrable limit \cite{Aubry}.
\end{rem}

\section{Proof of local connection theorems}

\label{sec:proof}

Following Shilnikov \cite{Shilnikov}, we will rewrite the boundary value problem (\ref{eq:BC0})  as  a fixed point problem.
First it is convenient to
make a change of variables.

\begin{prop}\label{prop:change}
There is a diffeomorphism $\Phi$ of  a neighborhood of $V\times (0,0)$ in $V\times\R^{k}\times \R^{k}$ such that:
\begin{itemize}
\item $\Phi$ is almost identity near $V$:
\begin{equation}
\label{eq:Phi}
\Phi(z,q,p)=(w,u,v)=(z,q,p)+O_2(q,p).
\end{equation}
\item
For any $z_0\in V_0\Subset V$,
\begin{equation}
\begin{array}{l}
\label{eq:PhiW}
 \Phi W^+_\loc(z_0)=\{(z_0,u,0):u\in B_r\},\\
\Phi W^-_\loc(z_0)=\{(z_0,0,v):  v\in B_r \}.  
\end{array}
\end{equation}
\item The flow
$\Phi\circ\phi^t\circ \Phi^{-1}$ on $\Phi W^\pm_\loc(z_0)$ is
linear:
\begin{equation}
\begin{array}{l}
\Phi\circ\phi^t\circ \Phi^{-1}(z_0,u,0)=(z_0,e^{-\lambda(z_0)t}u,0),\label{eq:flow}\\
\Phi\circ\phi^t\circ \Phi^{-1}(z_0,0,v)=(z_0,0,
e^{\lambda(z_0)t}v).
\end{array}
\end{equation}
\end{itemize}
\end{prop}

In general $\Phi$ is not symplectic.

\proof We modify the coordinates  $(z,q,p)$ in $U$ by setting
$$
u=q-f_-(z,p),\quad v=p-f_+(z,q),
$$
where $f_\pm$ are as in (\ref{eq:W+-}). In the variables $(z,u,v)$,
the local stable and unstable manifolds $W^\pm_\loc(V)$ are given by $v=0$
and $u=0$ respectively. Hence for $z_0\in V_0$, the manifold $W_\loc^+(z_0)$ is given by the
equations
$$
v=0,\quad z_0=z+\eta_+(z,u),
$$
and
$W_\loc^-(z_0)$ by the equations
$$
u=0,\quad z_0=z+\eta_+(z,v),
$$
where
$$
 \eta_+(z,u)=O_2(u),\quad \eta_-(z,v)=O_2(v).
$$
The projection $\pi_+$ is given by $z_0=z+\eta_+(z,u)$, and
the projection $\pi_-$ by $z_0=z+\eta_-(z,v)$.

We change the variable $z$ to
$$
w=z+\eta_+(z,u) + \eta_-(z,v).
$$
Then
\begin{eqnarray*}
w|_{W^+_\loc(z_0)}= z+\eta_+(z,u) + \eta_-(z,0) =z+\eta_+(z,u)=z_0,\\
w|_{W^-_\loc(z_0)}= z+\eta_+(z,0) + \eta_-(z,v) =z+\eta_-(z,v)=z_0.
\end{eqnarray*}
Thus the diffeomorphism $\Phi(z,q,p)=(w,u,v)$ satisfies
(\ref{eq:PhiW}).

The restriction of the Hamiltonian system to $W^+_\loc(w)$ is now
\begin{eqnarray}
\label{eq:restrict} \dot u = - \lambda(w)u + O_2(u) .
\end{eqnarray}
Since there are no resonances of order $\ge 2$, by Sternberg's theorem \cite{Sternberg},
there is a smooth normalizing transformation $u\to
\bar u =\phi(u,w)=u+O_2(u)$, smoothly depending on $w$ and
transforming system (\ref{eq:restrict}) to its
linear part $\dot {\bar u}  = - \lambda(w)\bar u$.    Similarly, we
can  transform  the system on $W_\loc^-(w)$ to $\dot {\bar v}  =
\lambda(w)\bar v$ via the change $v\to \bar v =\psi(v,w)=v+O_2(v)$. Then the map $\bar \Phi(z,q,p)=(w,\bar u,\bar v)$
satisfies (\ref{eq:flow}). We will
use the same notation $u,v$ for the new variables $\bar u,\bar v$.
Proposition \ref{prop:change} is proved.
\qed

\medskip

\begin{rem}
The last part of the proof is   the main place in the paper
where the equal eigenvalues case (\ref{eq:Apm}) differs from the generic case (\ref{eq:generic}).
Then there may be resonances, and the normal form is more complicated.
\end{rem}

The variables $u,v$ are closely related to the limit directions:
for $a_\pm=(z_\pm,q_\pm,p_\pm)\in W_\loc^\pm(z_0)$, we have
\begin{equation}
\label{eq:uv}
u(a_+)=v_+(z_0,q_+),\quad v(a_-)=v_-(z_0,p_-).
\end{equation}

In the variables $w,u,v$,   the Hamiltonian system   takes the
form
\begin{equation}
\label{eq:wuv}
\begin{array}{lcl}
\dot w &=& O(u;v),\\
\dot u &=& - \lambda(w)u + O(u;v),\\
\dot v&=&  \lambda(w)v + O(u;v).
\end{array}
\end{equation}

\begin{rem}
Here $O(u;v)$ means   a function of the form
$$
\sum_{|i|=|j|=1}a_{ij}(w,u,v)u^i v^j
$$
with smooth coefficients. Thus it vanishes on $W^+\cup W^-$.
\end{rem}

The Hamiltonian is transformed to
\begin{equation}
\label{eq:H(w,u,v)}
\calH(w,u,v)=H\circ \Phi^{-1}(w,u,v)=-\lambda(w)\langle u,v\rangle +O_3(u,v).
\end{equation}
However,  since $\Phi$ is
non-symplectic, system (\ref{eq:wuv}) does not have a  standard Hamiltonian form.

Finally we make a time change $d\tau=\lambda(w)\,dt$ and obtain
the system
\begin{equation}
\label{eq:final}
\begin{array}{lcl}
w' &=& O(u;v),\\
u' &=& - u + O(u;v),\\
v' &=&   v + O(u;v).
\end{array}
\end{equation}
Once a solution of system (\ref{eq:final})
is known, the time $t$ is determined by
\begin{equation}
\label{eq:t}
t=\theta(\tau)=\int_0^\tau\frac{ds}{\lambda(w(s))}.
\end{equation}

Next we reformulate Theorem \ref{thm:Shil} in the new variables.

\begin{prop} \label{prop:Shil-2}
Suppose $r>0$ is sufficiently small and $\calT\ge
1$.  Let $w_0\in V_0$ and $u_+,v_-\in B_r$. Then:
\begin{itemize}
\item
There exists a unique solution
\begin{equation}
\label{eq:uvw} \sigma(t)=(w(\tau),u(\tau),v(\tau))\in V\times B_{r}\times
B_{r},\qquad |\tau|\le \calT,
\end{equation}
of (\ref{eq:final}) satisfying the  initial-boundary conditions
\begin{equation}
\label{eq:bound3} w(0)=w_0,\quad u(- \calT)=u_+,\quad v(\calT)=v_-.
\end{equation}
\item
$\sigma$ smoothly depends on $(w_0,u_+,v_-,\calT)\in V_0\times B_r\times B_r\times[1,+\infty)$.
\item
Set
\begin{equation}
\label{eq:b+-}
(w_\pm,u_\pm,v_\pm)=\sigma(\mp \calT), \quad (w_0,u_0,v_0)=\sigma(0).
\end{equation}
As $\calT\to+\infty$, we have
\begin{equation}
\begin{array}{lcl}
u_0&=&u_+e^{-\calT} + e^{-2\calT}O(r^2),\label{eq:u(0)}\\
v_0&=&v_-e^{-\calT} + e^{-2\calT}O(r^2),\label{eq:v(0)}\\
w_+&=&w_0+\calT e^{-2\calT}O(r^2),\label{eq:w+}\\
w_-&=&w_0+\calT e^{-2\calT}O(r^2),\\
u_-&=&e^{-2\calT}(u_+ + O(r^2)),\label{eq:u-}\\
v_+&=&e^{-2\calT}(v_-+O(r^2)).\label{eq:v+}
\end{array}
\end{equation}
\item The initial and final time moments are
\begin{equation}
\label{eq:T+-}
T_\pm=\theta(\mp\calT)=\mp\lambda(w_0)\calT+\calT^2 e^{-2\calT}O(r^2).
\end{equation}
\end{itemize}
\end{prop}

\begin{rem}
The meaning of $O(r^2)$ is as in (\ref{eq:norm}): this is  a function $f(w_0,u_+,v_-,\calT)$
with $\|f\|_{C^1(V\times B_r\times B_r)}\le Cr^2$, where the constant is independent of $r$ and $\calT$,
and the weighted norm (\ref{eq:norm}) is used for the derivatives in $u_+$
and $v_-$.
\end{rem}

\medskip\noindent
{\it Proof of Proposition \ref{prop:Shil-2}.} We follow Shilnikov \cite{Shilnikov}. Set
\begin{equation}
\label{eq:u,v}
u=e^{-\tau-\calT}\xi, \quad v=e^{\tau-\calT}\eta.
\end{equation}
In the variables $w,\xi,\eta$, system (\ref{eq:final}) takes the form
\begin{equation}
\label{eq:final2}
\begin{array}{lclll}
w' &=&O(e^{-\tau-\calT}\xi;e^{\tau-\calT}\eta)&=& e^{-2\calT}O(\xi;\eta),\\
\xi' &=& e^{\tau+\calT}O(e^{-\tau-\calT}\xi;e^{\tau-\calT}\eta)&=&e^{\tau-\calT}O(\xi;\eta),\\
\eta' &=&
e^{\calT-\tau}O(e^{-\tau-\calT}\xi;e^{\tau-\calT}\eta)&=&e^{-\tau-\calT}
O(\xi;\eta).
\end{array}
\end{equation}
Here $O(\xi;\eta)$ is a function of the form
$$
\sum_{|i|=|j|=1}a_{ij}(w,\xi,\eta,\tau,\calT)\xi^i\eta^j,
$$
where the coefficients are smooth and uniformly bounded  for $\calT\ge 1$ and
$|\tau|\le\calT$.

Using  (\ref{eq:bound3}), we
obtain a system of integral equations
\begin{eqnarray}
\label{eq:int-w}
w(\tau)&=&w_0+ \int_0^\tau
e^{-2\calT}O(\xi(s);\eta(s))\,ds,\\
\label{e:int-xi}
\xi(\tau)&=&u_+ +\int_{-\calT}^\tau e^{s-\calT}O(\xi(s);\eta(s))\,ds,\\
\eta(\tau)&=&v_- +\int_{\calT}^\tau
e^{-s-\calT}O(\xi(s);\eta(s))\,ds. \label{eq:int-eta}
\end{eqnarray}

Let
$$
X=C^0([-\calT,\calT],\R^{2m}\times\R^{k}\times\R^{k}),
$$
be the Banach space with the norm
$$
\|(w,\xi,\eta)\|=\max\{\|w\|_{C^{0}},\|\xi\|_{C^0},\|\eta\|_{C^0}\},
$$
and let
$$
Y=\{(w,\xi,\eta):\|(w-w_0,\xi,\eta)\|\le 2r\}
$$
be a ball in $X$. We take $r>0$ so the small that the right hand sides of  equations
(\ref{eq:int-w})--(\ref{eq:int-eta}) are defined for
$(w,\xi,\eta)\in Y$. Then the right hand sides define a map $F:Y\to X$.

\begin{lem}
Let $r>0$ be sufficiently small. Then  $F(Y)\subset Y$ and
$F:Y\to Y$ is a contraction.
\end{lem}

\proof If $(w,\xi,\eta)\in Y$, then $|\xi(\tau)|,|\eta(\tau)|\le 2r$
for $|\tau|\le\calT$. There is a constant $C>0$, independent of $r$
and $\calT$, such that $|O(\xi;\eta)|\le Cr^2$. Set
$F(w,\xi,\eta)=(w_1,\xi_1,\eta_1)$. Then by
(\ref{eq:int-w})--(\ref{eq:int-eta}),
\begin{eqnarray*}
|\xi_1(\tau)-u_+|&\le&  \int_{-\calT}^\tau
e^{s-\calT}|O(\xi(s);\eta(s))|\,ds\le Cr^2(e^{\tau-\calT}-e^{-2\calT})\le r,\\
|\eta_1(\tau)-v_-|&\le&  \int_\tau^\calT
e^{-s-\calT}|O(\xi(s);\eta(s))|\,ds\le Cr^2(e^{-\tau-\calT}-e^{-2\calT})\le r,\\
|w_1(\tau)-w_0|&\le& e^{-2\calT}\bigg|\int_0^\tau Cr^2\,ds\bigg|=Cr^2\calT
e^{-2\calT}\le r,
\end{eqnarray*}
if $r<C^{-1}$. Hence $F(w,\xi,\eta)\in Y$. Similarly we show that for small $r>0$
the Lipschitz constant for $F$ is less than $1$, so $F$ is a contraction. \qed

\medskip

 Let $(w,\xi,\eta)\in Y$ be the
fixed point for $F$. Then by (\ref{eq:u,v}),
\begin{eqnarray}\label{eq:estimate1}
|u(\tau)-u_+e^{-\tau-\calT}|&=&e^{-\tau-\calT}|\xi(\tau)-u_+|\le Cr^2
e^{-2\calT}(1-e^{-\tau-\calT}),\\
|v(\tau)-v_-e^{\tau-\calT}|&=&e^{\tau-\calT}|\eta(\tau)-v_+|\le Cr^2
e^{-2\calT}(1-e^{-\tau-\calT}). \label{eq:estimate2}
\end{eqnarray}
We obtain
\begin{eqnarray*}
|w(\tau)-w_0|&\le  & Cr^2 \calT e^{-2\calT},\\
|u(\tau)-u_+e^{-\tau-\calT}|&\le & Cr^2e^{-2\calT},\\
|v(\tau)-v_-e^{\tau-\calT}|&\le & Cr^2e^{-2\calT}.
\end{eqnarray*}
Then by (\ref{eq:t}),
$$
|\theta(\tau)-\lambda(w_0)\tau| \le  Cr^2\calT^2 e^{-2\calT}.
$$

It remains to estimate the derivatives of solution $\sigma$ with
respect to $w_0,u_+,v_-$. Then we use  integral equations for the
corresponding variational system and get  e.g.\
$$
\bigg|\frac{\partial}{\partial u_+}u(\tau)-e^{-\tau-\calT}I\bigg|\le Cr
e^{-2\calT}.
$$
Similar estimates hold for other variables. This gives (\ref{eq:u(0)})
and then (\ref{eq:T+-}) follows from (\ref{eq:t}).

Last we check that $u(\tau),v(\tau)\in B_r$ for $|\tau|\le\calT$.
Equation (\ref{eq:estimate1}) gives
\begin{eqnarray*}
|u(\tau)|&\le & e^{-\tau-\calT}|u_+|+|u(\tau)-u_+e^{-\tau-\calT}|\\
&\le & r-(1-e^{-\tau-\calT})(r-Cr^2e^{-2\calT})\le r
\end{eqnarray*}
if $r<C^{-1}$. Thus $u(\tau)\in B_r$ for $|\tau|\le \calT$.
Similarly (\ref{eq:estimate2}) implies $v(\tau)\in B_r$ for  $|\tau|\le \calT$.

Proposition \ref{prop:Shil-2} is proved.
\qed

\medskip

Next we prove an analog of Theorem \ref{thm:Shil_mu} in the variables $w,u,v$.

 \begin{prop} \label{prop:uv_mu}
 Suppose $r>0$   and $\mu_0>0$ are sufficiently small.
 Let $w_0\in V_0$ and $(u_+,v_-)\in Q_r$. Then for every $\mu\in (0,\mu_0]$:
\begin{itemize}
\item
There exists $\calT>0$ 
and  a unique solution (\ref{eq:uvw}) with $\calH=\mu$ satisfying  (\ref{eq:bound3}).
\item
$\calT$ and $\sigma$ smoothly depend on $(w_0,u_+,v_-,\mu)\in V_0\times Q_r\times (0,\mu_0]$. Moreover
\begin{equation}
\label{eq:calT}
\calT=-\frac{1}{2\lambda(w_0)}\ln\left(-\frac{\mu}{\lambda(w_0)\langle
u_+,v_-\rangle}\right)+O(\sqrt{\mu})
\end{equation}
\item 
The boundary points (\ref{eq:b+-}) satisfy
\begin{equation}\label{eq:est4}
 \begin{array}{lcl}
w_+&=&w_0+O(\mu|\ln\mu|),\\
w_-&=&w_0+ O(\mu|\ln\mu|),\\
u_-&=&-\frac{\mu u_+}{\lambda(w_0)\langle u_+,v_-\rangle} + O(\mu),\\
v_+&=& -\frac{\mu v_-}{\lambda(w_0)\langle u_+,v_-\rangle}+O(\mu).
\end{array}
\end{equation}
\item
The initial and final time moments are
\begin{equation}
\label{eq:Tpm}
T_\pm=\theta(\mp\calT)= \mp\calT+O(\mu|\ln\mu|^2).
\end{equation}
\end{itemize}
\end{prop}

As before, $O(\mu)$ or $O(\mu|\ln\mu|)$ means a smooth function with  weighted $C^1$ norm bounded by
$C\mu$ or $C\mu|\ln\mu|$, where $C$ is independent of $r$ and $\mu$.

\proof
Equations (\ref{eq:H(w,u,v)}) and (\ref{eq:u(0)}) imply that on the solution (\ref{eq:uvw}),
$$
\calH|_\sigma=-e^{-2\calT}\lambda(w_0)\langle u_+,v_-\rangle +e^{-3\calT}O(r^3).
$$
For $\calH|_\sigma=\mu$, the implicit function theorem gives
$$
e^{-2\calT}=-\frac{\mu+O(\mu^{3/2})}{\lambda(w_0)\langle u_+,v_-\rangle},
$$
which implies (\ref{eq:calT}). Hence
$$
e^{-2\calT}O(r^2)=O(\mu),\quad \calT e^{-2\calT}O(r^2)=O(\mu|\ln\mu|).
$$
Then (\ref{eq:est4}) follow from  (\ref{eq:w+}), and (\ref{eq:Tpm}) from (\ref{eq:t}).
\qed

\medskip

\noindent{\it Proof of Theorem \ref{thm:xy_mu}.} We rewrite Proposition \ref{prop:uv_mu} in the variables $(x,y,q,p)$
via the change (\ref{eq:Phi}), where $z=(x,y)$.
Let
$$
x=X(w,u,v),\quad y=Y(w,u,v),\quad q=Q(w,u,v),\quad p=P(w,u,v)
$$
be the components of $\Phi^{-1}$.
According to (\ref{eq:est4}), to find a solution satisfying boundary conditions (\ref{eq:bound0}),
for given $x_+,y_-,q_+,p_-,\mu$ we need to find $w_0,u_+,v_-$ such that
\begin{eqnarray*}
X(w_0+O(\mu|\ln\mu|),u_+, -\frac{\mu v_+}{\lambda(w_0)\langle u_+,v_-\rangle} + O(\mu))=x_+,\\
Y(w_0+O(\mu|\ln\mu|),-\frac{\mu u_+}{\lambda(w_0)\langle u_+,v_-\rangle}+O(\mu) ,v_-)=y_-,\\
Q(w_0+O(\mu|\ln\mu|),u_+, -\frac{\mu v_-}{\lambda(w_0)\langle u_+,v_-\rangle} + O(\mu))=q_+,\\
P(w_0+O(\mu|\ln\mu|), -\frac{\mu u_+}{\lambda(w_0)\langle u_+,v_-\rangle} + O(\mu),v_-)=p_-.
\end{eqnarray*}
Hence
\begin{eqnarray}
X(w_0,u_+,0) + O(\mu|\ln\mu|)=x_+,\label{eq:X}\\
Y(w_0,0,v_-)+ O(\mu|\ln\mu|)=y_-,\label{eq:Y}\\
Q(w_0,u_+,0) + O(\mu)=q_+,\label{eq:Q}\\
P(w_0,0,v_-)+O(\mu)=p_-.\label{eq:P}
\end{eqnarray}
Equations (\ref{eq:Q})--(\ref{eq:P}) and (\ref{eq:u,v}) imply
\begin{eqnarray*}
u_+&=&v_+(w_0,q_+)+O(\mu),\\
v_-&=&v_-(w_0,p_-)+O(\mu).
\end{eqnarray*}
Then by (\ref{eq:X})--(\ref{eq:Y}),
$$
w_0=\zeta(x_+,y_-,q_+,p_-)+O(\mu|\ln\mu|).
$$
Let $\sigma(\tau)$ be the trajectory in Proposition \ref{prop:uv_mu} corresponding
to $w_0,u_+,v_-$ and let $t=\theta(\tau)$ be the corresponding time.  Set
$$
 t_0=\frac12(T_++T_-),\quad T=\frac12(T_--T_+).
$$
Then
\begin{equation}
\label{eq:back}
\gamma(t)=\Phi^{-1}(\sigma(\theta(t+t_0))),\qquad -T\le t\le T,
\end{equation}
satisfies the conditions of Theorem \ref{thm:xy_mu}.
\qed

\medskip

\noindent{\it Proof of Theorem
\ref{thm:x+y-}}.   Now we use Proposition \ref{prop:Shil-2}.
For given $x_+$, $y_-$, $q_+$, $p_-$, $T$  we need
to find $w_0$, $u_+$, $v_-$, $\calT$ such that
\begin{eqnarray*}
X(w_0+\calT e^{-2\calT}O(r^2),u_+,e^{-2\calT} v_-+e^{-2\calT} O(r^2))=x_+,\\
Y(w_0+\calT e^{-2\calT}O(r^2),e^{-2\calT} u_++e^{-2\calT} O(r^2),v_-)=y_-,\\
Q(w_0+\calT e^{-2\calT}O(r^2),u_+,e^{-2\calT} v_-+e^{-2\calT} O(r^2))=q_+,\\
P(w_0+\calT e^{-2\calT}O(r^2),e^{-2\calT} u_++e^{-2\calT} O(r^2),v_-)=p_-,\\
\lambda(w_0)\calT+\calT e^{-2\calT}O(r^2)=T.
\end{eqnarray*}
One can check that  for large $T$, this is possible by the implicit function
theorem. Let $\sigma(t)$ be the trajectory  (\ref{eq:uvw}).
Define $\gamma(t)$ as in (\ref{eq:back}). Theorem \ref{thm:x+y-} follows easily. \qed

\medskip

The proof of Theorem \ref{thm:Shil} is similar, and we skip it.

\section{Generating functions of the scattering map}

In this section we  relate the generating functions of the stable
and unstable manifolds $W^\pm$ and  of the
scattering map $\calF$.

Let $f:V^-\to V^+$ be a local branch of $\calF$ represented by
a generating function (\ref{eq:genf})
in symplectic coordinates $z_\pm=(x_\pm,y_\pm)$ in $V^\pm$.
Let  $(x_\pm,y_\pm,q_\pm,p_\pm)$ be the symplectic
coordinates in a tubular neighborhood
$$
U^\pm\cong V^\pm\times B_r\times B_r
$$
of $V^\pm$ such that the stable and unstable manifolds $W_\loc^\pm(V^\pm)$ are graphs (\ref{eq:W+-}).
As in (\ref{eq:Nmu}), take the cross sections
\begin{equation}
\begin{array}{l}
N^+=\{(z_+,q_+,p_+)\in U^+\cap\Sigma_0: q_+\in S_r\},\label{eq:N0}\\
N^-=\{(z_-,q_-,p_-)\in U^-\cap\Sigma_0: p_-\in S_r\}.
\end{array}
\end{equation}

Let  $\sigma$ be the transverse heteroclinic joining a point $c_0=(a_0,b_0)\in V^-$ with
$c_1=f(c_0)=(a_1,b_1)\in V^+$.
Let
$$
\sigma(t^\pm)=(a_\pm,b_\pm,c_\pm,d_\pm)=A^\pm\in N^\pm
$$
be the  intersection points  of $\sigma$ with $N^\pm$ such that
$\sigma(t)\in U^-$ for $t\le t^-$ and $\sigma(t)\in U^+$ for $t\ge
t^+$. Since $\sigma$ crosses $N^\pm$ transversely in $\Sigma_0$,
there exist neighborhoods $O^\pm$ of $A^\pm$  such that the
Poincar\'e map
\begin{equation}
\label{eq:calP}
\calP:O^-\cap N^-\to O^+\cap N^+,\qquad \calP(B)=\phi^{\tau(B)}(B),
\end{equation}
is a smooth symplectic  diffeomorphism. We have
$\tau(A^-)=t^+-t^-$ and $\calP(A^-)=A^+$.
We will locally represent $\calP$ by a generating function.

Suppose the neighborhoods $O^\pm$ are sufficiently small.
Let $D$ be a small neighborhood of $C=(b_-,d_-,a_+,c_+)$ and
$$
K=\{X=(y_-,p_-,x_+,q_+)\in D: q_+,p_-\in S_r\}.
$$

\begin{prop}\label{prop:conjugate}
The coordinates $x_+,q_+$ can be slightly modified in $O^+$ in
a way which does not invalidate the results of sections 3--4 and so that
\begin{itemize}
\item
For any $X=(x_-,q_-,y_+,p_+)\in K$ there exist points $B^\pm=(x_\pm,y_\pm,q_\pm,p_\pm)\in N^\pm\cap O^\pm$ such that
$\calP(B^-)= B^+$.
\item
$B^\pm=B^\pm(X)$  are smooth functions and $B^\pm(C)=A^\pm$.
\item
The Poincar\'e map (\ref{eq:calP}) is locally represented by a smooth
generating function $F(X)$ on $K$: for $B^\pm\in N^\pm\cap O^\pm$,
\begin{equation}
\label{eq:F0} \calP(B^-)=B^+\;\Leftrightarrow\; d
F(X)=p_+\, d q_+ + y_+\, d x_+ + x_-\, d y_- + q_-\,
d p_-.
\end{equation}
\end{itemize}
\end{prop}

\proof Consider the Lagrangian manifolds
\begin{eqnarray*}
L^+=\{(x_+,y_+,q_+,p_+)\in U^+: x_+=a_+,\; q_+=c_+\},\\
L^-=\{(x_-,y_-,q_-,p_-)\in U^-:  y_-=b_-,\; p_-=d_-\}.
\end{eqnarray*}
Since $d(H|_{L^\pm})(A^\pm)\ne 0$,
$\Pi^\pm=L^\pm\cap N^\pm\cap O^\pm$ are smooth Lagrangian manifolds in $N^\pm$.
We need to show that the Lagrangian manifolds  $\calP(\Pi^-)$ and $\Pi^+$
are transverse in $N^+$ at $A^+$, i.e.
\begin{equation}
\label{eq:transverse}
T_{A^+}\calP(\Pi^-)\cap T_{A^+}\Pi^+=\{0\}.
\end{equation}

Since $\vbf(A^+)$ is transverse to $N^+$, the symplectic space $T_{A^+}N^+$ is identified
with the quotient space $\calW=T_{A^+}\Sigma_0/\R\vbf(A^+)$. The Lagrangian subspace $T_{A^+}\Pi^+$ is identified with $\calV^+=(T_{A^+}L^+\cap T_{A^+}\Sigma_0)/\R\vbf(A^+)\subset \calW$, and $T_{A^+}\calP(\Pi^-)$ with a Lagrangian subspace $\calV^-\subset\calW$.

The transversality condition (\ref{eq:transverse}) is $\calV^-\cap\calV^+=\{0\}$.
This  can be achieved  by a slight perturbation of the manifold $L^+$
via local modification of the coordinates $x_+,q_+$ in a neighborhood of the point $A^+$.
Set
$$
\tilde x_+=x_++\frac{\partial}{\partial y_+}\phi(y_+,p_+),\quad \tilde q_+=q_++\frac{\partial}{\partial p_+}\phi(y_+,p_+),
$$
where $\phi$ is a small smooth function supported near $(b_+,d_+)$ such that $d\phi(b_+,d_+)=0$.
Let   $\calH$ be the Hessian matrix of $\phi$ at $(b_+,d_+)$.
If we use the coordinates $\tilde x_+,y_+,\tilde q_+,p_+$, the manifold $L^+$ is replaced
by
$$
L_\phi^+=\{(x_+,y_+,q_+,p_+): x_++\frac{\partial}{\partial y_+}\phi(y_+,p_+)=a_+,\; q_++\frac{\partial}{\partial p_+}\phi(y_+,p_+)=c_+\}.
$$
Then $T_{A^+}L^+$ is replaced by a Lagrangian subspace $\calL^+_\calH=T_{A^+}L_\phi^+$
depending on $\calH$. Changing $\calH$, we get an open set $\{\calL_H^+\}$ of
Lagrangian subspaces in $T_{A^+}\Sigma_0$.
Hence we obtain an open set $\{\calV_{\calH}^+\}$ of Lagrangian subspaces $\calV_{\calH}^+=(\calL^+_\calH\cap T_{A^+}\Sigma_0)/\R\vbf(A^+)$ in $\calW$.
Thus for almost all $\calH$, the Lagrangian subspaces $\calV_{\calH}^+$ and $\calV^-$ are transverse.
\qed

\medskip

Proposition \ref{prop:conjugate} is more clear in local symplectic
coordinates $x_\pm,y_\pm,\xi_\pm,\eta_\pm$ on $N^\pm\cap O^\pm$ defined as in
(\ref{eq:local}). Then
\begin{equation}
\label{eq:Bpm}
 B^\pm\leftrightarrow (x_\pm,y_\pm,\xi_\pm,\eta_\pm),\quad X\leftrightarrow(y_-,\eta_-,x_+,\xi_+).
 \end{equation}
 Let
\begin{equation}
\label{eq:xi-eta}
x_+=x_+(x_-,y_-,\xi_-,\eta_-),\quad
\xi_+=\xi_+(x_-,y_-,\xi_-,\eta_-).
\end{equation}
be the components of the Poincar\'e map
\begin{equation}
\label{eq:Pcoord}
(x_-,y_-,\xi_-,\eta_-)\to (x_+,y_+,\xi_+,\eta_+).
\end{equation}
 Then the transversality condition (\ref{eq:transverse}) is
\begin{equation}
\label{eq:det}
\det \frac{\partial(x_+,\xi_+)}{\partial( x_-,\xi_-)}\bigg|_{A^-}
\ne 0.
\end{equation}
Under condition (\ref{eq:det}), equations (\ref{eq:xi-eta}) can be solved for
$$
x_-=x_-(y_-,\eta_-,x_+,\xi_+),\quad \xi_-=\xi_-(y_-,\eta_-,x_+,\xi_+),
$$
which gives the point $B^-(X)$ and then $B^+(X)=\calP(B^-)$.
The Poincar\'e map (\ref{eq:Pcoord})
is represented by the generating function
$$
\varphi(y_-,\eta_-,x_+,\xi_+)=F(y_-,p_-(\eta_-),x_+,q_+(\xi_+))
$$
as follows:
\begin{equation}
\label{eq:barF0} \calP(B^-)=B^+\;\Leftrightarrow\; d \varphi=\eta_+\,
d \xi_+ + y_+\, d x_+ + x_-\, d y_- + \xi_-\, d \eta_-.
\end{equation}

\begin{rem}
Transversality of $\sigma$ implies, without any modification of the coordinates,
that $\calP$ can be represented by a generating function of the variables $x_-,q_-,y_+,p_+$.
However,  for the proof of Theorem \ref{thm:period}  the generating function of the variables $y_-,p_-, x_+,q_+$
is more convenient.
\end{rem}

Let $S_\pm$ be the   generating functions (\ref{eq:Spm}) of the local stable and unstable manifolds $W^\pm$.
Set
\begin{equation}
\label{eq:G}
G_{x_0,y_1}(X)=S_-(x_0,y_-,p_-)-
F( X)+S_+(x_+,y_1,q_+),\qquad X\in K.
\end{equation}

\begin{prop}
\label{prop:nondeg}
\begin{itemize}
\item
$X\in K$ is a critical point of $G_{x_0,y_1}$ iff
$B^-(X)\in W^-(x=x_0)$ and $B^+(X)\in W^+(y=y_1)$, i.e.\ the points $B^\pm$ lie on a heteroclinic
orbit.
\item
If the  heteroclinic orbit $\sigma$  is transverse, then $C$ is a
nondegenerate critical point of $G_{a_0,b_1}$ on
$K$.
\item
For $(x_0,y_1)$ close to $(a_0,b_1)$, the function
$G_{x_0,y_1}$ has a nondegenerate critical point $X(x_0,y_1)\in K$ such that $X(a_0,b_1)=C$. The
critical value is the generating function of the scattering map:
\begin{equation}
\label{eq:G2}
S(x_0,y_1)=\Crit_{X\in K}\,G_{x_0,y_1}(X)=G_{x_0,y_1}(X(x_0,y_1)).
\end{equation}
\end{itemize}
\end{prop}

 \proof
 We represent $X\in K$ and the corresponding points $B^\pm(X)$ in local coordinates 
 as in (\ref{eq:Bpm}).  Set
 \begin{eqnarray*}
 &R_-(x_0,y_-,\eta_-)=S_-(x_0,y_-,p_-(\eta_-)),\\
 &R_+(x_+,y_1,\xi_+)=S_+(x_+,y_1,q_+(\xi_+)),\\
 &R(x_0,y_1,y_-,\eta_-,x_+,\xi_+)=G_{x_0,y_1}(y_-,p_-(\eta_-),x_+,q_+(\xi_+)).
 \end{eqnarray*}
 Then by (\ref{eq:barF0}),
 \begin{eqnarray}
 &dR=(\hat x_--x_-)\,dy_-+(\hat \xi_--\xi_-)\,d\eta_-
 +(\hat y_- - y_-)\,dx_++(\hat \eta_+-\eta_-)\,d\xi_+ \nonumber\\
 &+y_0\,dx_0+x_1\,dy_1,\label{eq:hat}
 \end{eqnarray}
 where
 \begin{eqnarray*}
 \hat x_-=\frac{\partial}{\partial y_-}R_-(x_0,y_-,\eta_-),\quad \hat \xi_-=\frac{\partial}{\partial\eta_-}R_-(x_0,y_-,\eta_-),\\
 \hat y_+= \frac{\partial}{\partial x_+}R_+(y_1,x_+,\xi_+),
 \quad \hat \eta_+=\frac{\partial}{\partial \xi_+}R_+(y_1,x_+,\xi_+).
 \end{eqnarray*}
 Let
 \begin{equation}
 \label{eq:hatB}
 \hat B^-\leftrightarrow(\hat x_-,y_-,\hat \xi_-,\eta_-),\quad  \hat B^+\leftrightarrow(x_+,\hat y_+,\xi_+,\hat \eta_+).
 \end{equation}
 By (\ref{eq:pi+})--(\ref{eq:pi-}), $\hat B^-\in W^-(x=x_0)$ and $\hat B^+\in W^+(y=y_1)$.
 If $X$ is a critical point of $G_{x_0,y_1}$,
 then $\hat B^\pm=B^\pm$.
 Hence $B^\pm$ lie on a heteroclinic orbit which proves the first item of Proposition \ref{prop:nondeg}.
 Then by (\ref{eq:hat}),
 \begin{equation}
 \label{eq:dG}
 dR=y_0\,dx_0+x_1\,dy_1.
 \end{equation}

 Suppose that $C$ is a degenerate critical point of   $G_{a_0,b_1}$ on $K$.
 Then there is a family of nearly critical points
 $$
 X(\eps)\leftrightarrow(y_-(\eps),\eta_-(\eps),x_+(\eps),\xi_+(\eps))
 $$
 such that $X(0)=C$, $X'(0)\ne 0$  and
 \begin{equation}
 \label{eq:eps}
 dG_{a_0,b_1}(X(\eps))=O(\eps^2).
 \end{equation}
 Let
$$
B^\pm(\eps)\leftrightarrow(x_\pm(\eps),y_\pm(\eps),\xi_\pm(\eps),\eta_\pm(\eps))
$$
be the points corresponding to $X(\eps)$ by Proposition  \ref{prop:conjugate}
and let
$$
 \hat B^-(\eps)\leftrightarrow(\hat x_-(\eps),y_-(\eps),\hat \xi_-(\eps),\eta_-(\eps)),\quad
 \hat B^+(\eps)\leftrightarrow(x_+(\eps),\hat y_+(\eps),\xi_+(\eps),\hat \eta_+(\eps))
 $$
be the points defined in (\ref{eq:hatB}).  Then (\ref{eq:hat}) and (\ref{eq:eps}) imply $\hat B^\pm(\eps)=B^\pm(\eps)+O(\eps^2)$.
Applying the Poincar\'e map, we obtain
$$
\calP(\hat B^-(\eps))=\calP(B^-(\eps))+O(\eps^2)=B^+(\eps)+O(\eps^2)\in
W^-(x=a_0).
$$
Thus  the curve $B^+(\eps)\in W^+(y=b_1)$ is tangent to  $W^-(x=a_0)$.
This contradicts the assumption that the heteroclinic $\sigma$ is transverse.

The last item follows from the first two and (\ref{eq:dG}).
 \qed

\medskip

Suppose now that $\mu_0>0$ is sufficiently small and let
$\mu\in[-\mu_0,\mu_0]$. We introduce cross sections
$N_\mu^\pm\subset\Sigma_\mu\cap U^\pm$ as in (\ref{eq:Nmu}). Then $N_0^\pm=N^\pm$.
By the implicit function theorem,
the Poincar\'e map $\calP_\mu:O^-\cap N_\mu^-\to O^+\cap N_\mu^+$ is
well defined and coincides with $\calP$ for $\mu=0$. Proposition \ref{prop:conjugate} implies

\begin{cor}\label{cor:conjugate}
For any $\mu\in[-\mu_0,\mu_0]$ and $X=(y_-,p_-,x_+,q_+)\in K$:
\begin{itemize}
\item There exist $x_-,p_-,y_+,q_+$ such that $B^\pm(X,\mu)=(x_\pm,y_\pm,q_\pm,p_\pm)\in\Sigma_\mu$
and  $P_\mu(B^-)= B^+$.
\item
The Poincar\'e map $\calP_\mu:N_\mu^-\cap O^-\to N_\mu^+\cap O^+$  has a smooth
generating function $F_\mu(X)=F(X)+O(\mu)$, $X\in K$, smoothly
depending on $\mu\in[-\mu_0,\mu_0]$:
$$
\calP_\mu(B^-)=B^+\quad\Leftrightarrow\quad  d F_\mu(X)=p_+\, d q_+ + y_+\, d x_+ + x_-\, d
y_-+q_-\, d p_-.
$$
\end{itemize}
\end{cor}

\section{Variational problem}

In this section we define 2 functionals: one whose critical points correspond to
periodic heteroclinic chains and another whose critical points correspond to shadowing orbits on $\Sigma_\mu$.
Then Theorem \ref{thm:period} follows easily.

Let $c_{i+1}=f_i(c_i)$   be a $n$-periodic orbit of $\calF$ and let $\sigma=(\sigma_i)$ be the corresponding periodic heteroclinic chain: $c_i=\sigma_i(-\infty)$ and $c_{i+1}=\sigma_{i}(+\infty)$. In the symplectic coordinates $z_i=(x_i,y_i)$  in a
neighborhood $V_i$ of $c_i=(a_i,b_i)$,  $f_i$ is represented by a generating function
$S_i(x_i,y_{i+1})$ as in (\ref{eq:genfj}). Then
$\cbf=(c_i)_{i=0}^{n-1}$ is a  critical point of
the action functional (\ref{eq:calA}).

In a neighborhood $U_i\cong V_i\times B_r\times B_r$ of $c_i$ in
$ \calM $ we will use symplectic coordinates $(x_i,y_i,q_i,p_i)$
as in (\ref{eq:W+-}). Define the cross sections as in (\ref{eq:N0}):
\begin{eqnarray*}
N_i^+=\{(x_i,y_i,q_i,p_i)\in U_i\cap\Sigma_0:q_i\in S_r\},\\
N_i^-=\{(x_i,y_i,q_i,p_i)\in U_i\cap\Sigma_0:p_i\in S_r\}.
\end{eqnarray*}
Let
\begin{eqnarray*}
A_i^-=(a_i^-,b_i^-,c_i^-,d_i^-)\in N_i^-,\quad
A_{i+1}^+=(a_{i+1}^+,b_{i+1}^+,c_{i+1}^+,d_{i+1}^+)\in N_{i+1}^+
\end{eqnarray*}
be the first and last  intersection points of $\sigma_i$ with   $N_i^-$
and $N_{i+1}^+$ respectively.
Take small neighborhoods $O_i^\pm$  of $A_i^\pm$ and let
$\calP_i:N_i^-\cap O_i^-\to N_{i+1}^+\cap O_{i+1}^+$ be the local Poincar\'e map.
Then $\calP_i(A_i^-)=A_{i+1}^+$.

Let $D_i$ be a small neighborhood of $C_i=(b_i^-,d_i^-,a_{i+1}^+,c_{i+1}^+)$
and 
$$
K_i=\{X_i=(y_i^-,p_i^-,x_{i+1}^+,q_{i+1}^+)\in D_i:p_i^-, q_{i+1}^+\in S_r\}.
$$
By Proposition \ref{prop:conjugate}, without loss of generality we
may assume that for any $X_i=(y_i^-,p_i^-,x_{i+1}^+,q_{i+1}^+)\in
K_i$ there exist $x_i^-$, $q_i^-$, $y_{i+1}^+$, $p_{i+1}^+$, smoothly depending on $X_i$, such
that the points
\begin{eqnarray*}
B_i^-(X_i)=(x_i^-,y_i^-,q_i^-,p_i^-)\in N_i^-,\quad
B_{i+1}^+(X_i)=(x_{i+1}^+,y_{i+1}^+,q_{i+1}^+,p_{i+1}^+)\in N_{i+1}^+
\end{eqnarray*}
satisfy $\calP_i(B_i^-)=B_{i+1}^+$.
The Poincar\'e map $\calP_i$ is locally given by the generating function
$F_i(X_i)$ on $K_i$:
$$
 d F_i(X_i)=p_{i+1}^+\, d q_{i+1}^+ + y_{i+1}^+\, d x_{i+1}^+
+ x_i^-\, d y_i^-+ q_i^-\, d p_i^-.
$$

As in (\ref{eq:G}), let
\begin{eqnarray*}
G_i(x_i,y_{i+1},X_i )=S_i^-(x_i,y_i^-,p_i^-)- F_i(X_i)
+  S_{i+1}^+(x_{i+1}^+,y_{i+1},q_{i+1}^+).
\end{eqnarray*}
By Proposition \ref{prop:nondeg}, $X_i\to G_i(x_i,y_{i+1},X_i )$ has a nondegenerate critical value
\begin{equation}
\label{eq:Crit2}
S_i(x_i,y_{i+1})=\Crit_{X_i\in K_i}G_i(x_i,y_{i+1},X_i )= G_i(x_i,y_{i+1},X_i(x_i,y_{i+1} ))
\end{equation}
which is the generating
function of the symplectic map $f_{i}$.

Let
\begin{eqnarray*}
\calB(\zbf,\Xbf) =\sum_{i=0}^{n-1}( G_i(x_i,y_{i+1},X_i)-\langle
x_i,y_i\rangle),\qquad \zbf=(z_i)_{i=0}^{n-1},\quad \Xbf= (X_i)_{i=0}^{n-1},
\end{eqnarray*}
where
$$
z_i=(x_i,y_i)\in V_i,\quad X_i=(y_i^-,p_i^-,x_{i+1}^+,q_{i+1}^+)\in K_i,
$$
and 
$$
y_n=y_0,\quad x_n^+=x_0^+,\quad  q_n^+=q_0^+.
$$
In fact $\calB$ is a modified Maupertuis action of the concatenation of trajectories
of the Hamiltonian system on $\Sigma_0$. It is a smooth function on
$$
\calN=\calV\times\calK,\qquad
\calV=\prod_{i=0}^{n-1}V_i,\quad
\calK=\prod_{i=0}^{n-1}K_i.
$$

\begin{prop}
\label{prop:crit1}
\begin{itemize}
\item
For any $\zbf\in\calV$ close to $\cbf$, the function $\Xbf\in\calK\to\calB(\zbf,\Xbf)$
 has a nondegenerate critical point $\Xbf(\zbf)$.
The critical value  equals the action functional (\ref{eq:calA}):
$$
\calA(\zbf)=\Crit_{\Xbf\in\calK}\calB(\zbf,\Xbf)=\calB(\zbf,\Xbf(\zbf)).
$$
\item
Let
$(\cbf,\Cbf)$, $\Cbf=\Xbf(\cbf)$,  be the  critical point of $\calB$ corresponding to
the periodic orbit $\cbf$. If $\cbf$ is
nondegenerate, then $(\cbf,\Cbf)$ is nondegenerate.
\end{itemize}
\end{prop}

The first statement follows from (\ref{eq:Crit2}), and the second
from the following elementary and well known

\begin{lem}\label{lem:crit}
Let $f(x,y)$ be a smooth function and let let $y=h(x)$ be a nondegenerate critical point of $f(x,y)$
with respect to $y$.
Then $(x_0,y_0)$ is a nondegenerate critical point of $f(x,y)$ iff $x_0$ is a nondegenerate
critical point of $g(x)=f(x,h(x))$.
\end{lem}

Suppose now that the heteroclinic chain $\sigma$ is positive.
Let $\kappa>0$ and $r>0$ be so small that
$$
\langle c_i^+,d_i^-\rangle < -\kappa r^2.
$$
Then $(q_i^+,p_i^-)\in Q_r$ for $(q_i^+,p_i^-)$ close to $(c_i^+,d_i^-)$.

Take small $\mu_0>0$ and let $\mu\in(0,\mu_0]$.
Let  $R_i^\mu(Z_i)$, $Z_i=(x_i^+,y_i^-,q_i^+,p_i^-)\in V_i\times Q_r$, be the generating function in Theorem \ref{thm:xy_mu}
corresponding to $V_i\subset M$. It generates the Poincar\'e map $P_i^\mu: N_{i,\mu}^+\cap O_i^+\to N_{i,\mu}^-\cap O_i^-$  of the cross sections  $N_{i,\mu}^\pm\subset U_i\cap \Sigma_\mu$ defined in (\ref{eq:Nmu}).

Let $F_i^\mu(X_i)$,
$X_i=(y_i^-,p_i^-,x_{i+1}^+,q_{i+1}^+)\in K_i$,
be the generating function of the  Poincar\'e map $\calP_i^\mu:N_{i,\mu}^-\to N_{i+1,\mu}^+$ in Corollary
\ref{cor:conjugate}.
Set
$$
\calA_\mu(\Xbf)=\sum_{i=0}^{n-1}
(F_i^\mu(X_i)+R_i^\mu(Z_i)),\qquad \Xbf=(X_i)_{i=0}^{n-1}.
$$
We obtain

\begin{prop}
$\Xbf$ is a critical point of $\calA_\mu$ iff the corresponding points $B_i^\pm=B_i^\pm(X_i,\mu)\in N_{i,\mu}^\pm$ in Corollary
\ref{cor:conjugate} lie on a periodic orbit $\gamma_\mu$ in $\Sigma_\mu$.  Equivalently, $B_0^-$ is a fixed point
of the total Poincar\'e map
$$
P_{n-1}^\mu\circ \calP_{n-2}^\mu\circ\cdots\circ\calP_1^\mu\circ P_1^\mu\circ\calP^\mu:N_{0,\mu}^-\to N_{0,\mu}^-.
$$
\end{prop}

For $\mu=0$ we have
$$
\calA_0(\Xbf)=\sum_{i=0}^{n-1}
(F_i(X_i)+L_i(Z_i)),
$$
where $F_i(X_i)$ is the generating function of the Poincar\'e map $P_i$, and 
$L_i(Z_i)$ the generating function of the symplectic relation in (\ref{eq:L}):
\begin{equation}
\label{eq:L_i}  dL_i(Z_i)=y_i^+ d x_i^++y_i^-d x_i^-+p_i^+ d
q_i^+ + q_i^- d p_i^-.
\end{equation}
Proposition \ref{prop:reflect} implies that to  $\Xbf\in\calK$ there corresponds $\zbf(\Xbf)\in\calV$ such that
$$
\calA_0(\Xbf)=\Crit_{\zbf}\calB(\zbf,\Xbf)=\calB(\zbf(\Xbf),\Xbf).
$$
By Lemma \ref{lem:crit}, if $(\cbf,\Cbf)$ is a nondegenerate critical point of $\calB$ on $\calN$,
then $\cbf$ is a nondegenerate critical point of $\calA$, and $\Cbf$ is a nondegenerate critical point of $\calA_0$.

Now we can prove Theorem \ref{thm:period}.
Let $\cbf$ be a nondegenerate  periodic orbit of $\calF$ corresponding to a positive
heteroclinic chain $\sigma$.
By Proposition \ref{prop:crit1} it defines a nondegenerate critical point $(\cbf,\Cbf)$
of $\calB$ which gives a nondegenerate critical point $\Cbf$ of $\calA_0$.
By Theorem \ref{thm:Shil_mu},
\begin{equation}
\label{eq:calB}
\|\calA_\mu-\calA_0\|_{C^2}\le \const \, \mu|\ln\mu|.
\end{equation}
Hence for small $\mu>0$, $\calA_\mu$ has  a nondegenerate critical point $\Cbf_\mu=\Cbf+O(\mu|\ln\mu|)$
which gives a periodic shadowing trajectory $\gamma_\mu$. Theorem \ref{thm:period} is proved.
\qed

\begin{rem}
The constant in (\ref{eq:calB}) may depend on $n$, so in this proof we are unable to pass to the limit as $n\to+\infty$.
To get chaotic shadowing trajectories and prove Theorem \ref{thm:skew}, we need to use the $L_\infty$ norm on the space
of sequences. This will be done in a subsequent publication.
\end{rem}

\end{document}